\documentclass[a4paper,11pt, reqno]{amsart} 
\usepackage{amssymb,amsthm,amsmath}

\usepackage{amsfonts}
\usepackage{amscd}
\usepackage{amssymb}
\usepackage{enumerate}
\usepackage{graphicx}
\allowdisplaybreaks
\usepackage{mathabx}
\usepackage{color}
\usepackage{amsbsy}
\usepackage{graphicx}
\usepackage{amsthm}
\usepackage{amsmath}
\usepackage{amsxtra}
\usepackage{mathrsfs}
\usepackage{bbm}
\usepackage{dsfont}

\usepackage{ifthen}

\usepackage[colorlinks,citecolor=red,pagebackref,hypertexnames=false]{hyperref} 
\newtheorem{theorem}{Theorem}[section]
\newtheorem{cor}[theorem]{Corollary}
\newtheorem{lem}[theorem]{Lemma}

\newtheorem{prop}[theorem]{Proposition}
\newtheorem{prob}[theorem]{Problem}
\theoremstyle{definition}
\newtheorem{defn}[theorem]{Definitions}
\theoremstyle{remark}
\newtheorem{rem}[theorem]{Remark}
\numberwithin{equation}{section}
\definecolor{red}{rgb}{1.0, 0.0, 0.0}
\newcommand{\Bea}{\begin{eqnarray*}}
\newcommand{\Eea}{\end{eqnarray*}}
\newcommand{\Be} {\begin{equation*}}
\newcommand{\Ee} {\end{equation*}}
\newcommand{\be} {\begin{equation}}
\newcommand{\ee} {\end{equation}}
\newcommand{\bea} {\begin{eqnarray}}
\newcommand{\eea} {\end{eqnarray}}

\newcommand{\R}{\mathbb R}

\newcommand{\si}{\sigma}

\newcommand{\rc}{\mathfrak{R}}
\newcommand{\Rc}{\mathcal{R}}

\newcommand{\ecs}{\mathcal{E}_S}
\newcommand{\ecso}{\mathcal{E}_{S_1}}
\newcommand{\ecstw}{\mathcal{E}_{S_2}}
\newcommand{\ecsth}{\mathcal{E}_{S_3}}
\newcommand{\Sc}{\mathcal S}

\newcommand{\la}{\lambda}
\newcommand{\C}{{\mathbb C}}

\newcommand{\La}{\langle}
\newcommand{\Ra}{\rangle}

\newcommand{\Hc}{\mathcal{H}}

\newcommand{\om}{ \omega}
\newcommand{\D}{\Delta}
\newcommand{\Dk}{\Delta_{\kappa}}

\newcommand{\ap}{\alpha}
\newcommand{\bt}{\beta }

\newcommand{\K}{\kappa }
\newcommand{\gm}{\gamma}
\newcommand{\z}{\zeta}
\newcommand{\Gm}{\Gamma}
\newcommand{\ve}{\varepsilon}
\newcommand{\gk}{\gamma_\kappa}

\newcommand{\s}{\mathfrak{S}}

\renewcommand{\Re}{\operatorname{Re}}
\renewcommand{\Im}{\operatorname{Im}}

\DeclareMathOperator{\sgn}{sgn}

\title[Restriction Theorems For Fourier-Dunkl Transform]
{Restriction theorem for Fourier-Dunkl transform II: Paraboloid, sphere, and hyperboloid surfaces}
\author[P Jitendra K. Senapati, Pradeep B, S. S. Mondal and H. Mejjaoli]{P Jitendra Kumar Senapati, Pradeep Boggarapu, 
Shyam Swarup Mondal and Hatem Mejjaoli}
\address[P Jitendra K. Senapati.]{Department of Mathematics\\
   BITS Pilani K K Birla Goa Campus\\
    Zuarinagar, South Goa\\
403 726, Goa, India}
\email{jitusnpt@gmail.com}

\address[Pradeep B.]{Department of Mathematics\\
   BITS Pilani K K Birla Goa Campus\\
    Zuarinagar, South Goa\\
403 726, Goa, India}
\email{pradeepb@goa.bits-pilani.ac.in}

\address[S. S. Mondal]{Department of Mathematics\\
  Indian Institute Of Technology Delhi\\
    Delhi 11016, India}
\email{mondalshyam055@gmail.com}
\address[H. Mejjaoli]{Department of Mathematics\\
  Taibah University, College of Sciences,  PO BOX 30002 Al Madinah AL Munawarah, Saudi Arabia}
\email{hmejjaoli@gmail.com}

\begin{document}

\maketitle
%--------------------------------------------------------------------------------------------------------------------------------------------------
\begin{abstract}
This is a continuation of the paper \cite{JPSH} where the authors  introduced  and studied  the Fourier-Dunkl transform  on $\mathbb{R}^{n}\times \mathbb{R}^{d}$.  The main novelty of this paper is that  we here  prove  Strichartz's restriction theorem for the Fourier-Dunkl transform for certain   surfaces,  namely, paraboloid, sphere, and hyperboloid and its generalisation to the family of orthonormal functions.  Finally, as an application of these    restriction theorems, we establish   versions of  Strichartz estimates for orthonormal families of initial data  associated with Schr\"odinger's propagator  in the case of the Dunkl Laplacian and Klein-Gordon operator.

% We establish new Strichartz estimates for orthonormal families of initial data in the case of
% versions of Strichartz estimates
 
 %We establish new Strichartz estimates for orthonormal families of initial data in the case of the wave, Klein–Gordonand fractional Schrödinger equations

%This is a continuation of the paper \cite{JPSH}.  In this article, we prove Strichartz's restriction theorem for the Fourier-Dunkl transform which is defined by Jitendra K. S. et. al in \cite{JPSH} for the certain surfaces namely, paraboloid, sphere, and hyperboloid and its generalisation to the family of orthonormal functions. As an application of these restriction theorems, we derive the Strichartz inequalities associated with Schr\"odinger's propagator of Dunkl Laplacian and Klein-Gordon operator for the family of orthonormal functions.
\end{abstract}
%----------------------------------------------------------------------------------------------------------------------------------------------------
	
	\tableofcontents

\section{Introduction}
One of the most  important and  classical topic in harmonic analysis is the  restriction problem.  Let    $S\subset\R^n$ be a hyper-surface endowed with Lebesgue measure $d\mu,$ then one may ask the following:

\begin{prob}\label{problem1}
 For which exponents $1 \leq p \leq 2$, the Fourier transform of a function $f \in L^{p}\left(\mathbb{R}^{n}\right)$ belongs to $L^{2}(S)$, where $S$ is endowed with its $(n-1)$-dimensional Lebesgue measure $d \sigma$? More preciously,  for which exponents $1 \leq p \leq 2$, we have the inequality
\begin{equation}\label{Restriction1}
	\int_S|\mathcal{F}f(\xi)|^2 d\mu(\xi)\leq C \|f\|^2_{L^p(\R^n)}.
\end{equation}
\end{prob}
Here  $\mathcal{F}f$ denotes the Fourier transform for a   Schwartz class function $f$ on $\R^n$ and  is defined as
$$\mathcal{F}f(\xi)=\frac{1}{(2\pi)^{n/2}}\int_{\R^n}f(x)e^{-ix\cdot\xi}dx.$$
Note that Problem \ref{problem1} is famously known as   
 the restriction theorem for the Fourier transform.

% The phenomenon of restriction of the Fourier transform has attracted considerable interest in the thirty-five years since it was first discovered. Restriction theorems highlight the delicate connection between curvature and Fourier analysis, melding geometry and cancellation to great effect. In addition, there are many practical applications of restriction theorems to areas of harmonic analysis (Bochner-Riesz summability, for example) and PDE.
 
 The phenomenon of restriction of the Fourier transform attracted considerable interest  since it was first discovered.    E. M. Stein  first introduced the restriction phenomena    for studying the boundedness of the Fourier transform of an $L^p$-function in the Euclidean space $\mathbb{R}^n$ for some $n  \geq2$.      The restriction problem has a delicate connection to many other conjectures, notably the Kakeya and Bochner-Riesz conjectures. Furthermore, it  also closely related to that of estimating solutions to linear PDE such as the wave and Schr\"odinger equations;  and this connection was first observed by Strichartz in \cite{RS}. More specifically,    Strichartz inequalities for the Schr\"odinger and wave equations correspond to Fourier restriction estimates on the paraboloid and the cone, respectively. 
 
% It plays a vital role in the theory of partial differential equations and harmonic analysis.

  %z [76], leading to the family of estimates which bear his name.  Furthermore, t 

%The restriction theorem for the Fourier transform is a classical problem in harmonic analysis.   The restriction is deep and far-reaching and has a surprising connection to many other important problems in mathematical analysis, geometric measure theory,  number theory, combinatorics, harmonic analysis,   including several  conjectures such as the Bochner-Riesz conjecture, Kakeya conjecture,     as well as PDE conjectures such as the local smoothing conjecture, see \cite{tao}.

%
% However, we first recall the simplest case. For a given Schwartz class function $f$ on $\R^n,$ its Fourier transform is defined as
%$$\mathcal{F}f(\xi)=\frac{1}{(2\pi)^{n/2}}\int_{\R^n}f(x)e^{-ix\cdot\xi}dx.$$
%% \textcolor{red}{It is well know that $\mathcal{F}f$ is also a Schwartz class function, satisfying $\|f\|_{L^2(\R^n)}=\|\mathcal{F}f\|_{L^2(\R^n)}$ and $f=\mathcal{F}^2f(-x).$ } 
%Let $S\subset\R^n$ be a hyper-surface endowed with Lebesgue measure $d\mu,$ then one may ask the following (Restriction problem):
%
%\noindent{\bf \textit{Problem 1.}} For $n\geq 2,$ for which values of $1\leq p\leq 2,$ is it true that for all $f\in\Sc(\R^n)$ $$\int_S|\mathcal{F}f(\xi)|^2 d\mu(\xi)\leq C \|f\|^2_{L^p(\R^d)}?$$
Let the operator $\rc_S$ defined as $\rc_Sf=\mathcal{F}f|_S$, known as restriction operator. Then the Problem \ref{problem1} is equivalent to the boundedness of $\rc_S: L^p(\R^n) \to L^2(S).$ The operator $\ecs$ dual to the restriction operator is defined as
$$\ecs f(x)=\frac{1}{(2\pi)^{n/2}}\int_S f(\xi)e^{i\xi\cdot x} d\mu(\xi), \quad x\in\R^n,$$ for $f\in L^1(S).$
When  q is the conjugate exponent of  $p$, i.e., $\frac{1}{p}+\frac{1}{q}=1,$ by duality, Problem \ref{problem1}   can also be posed as follows:
 \begin{prob}\label{problem2}
 	For $n\geq 2,$ for which exponent $2\leq q\leq \infty,$ is it true that for all $f\in L^2(S)$, we have  $$\|\ecs f\|_{L^q(\R^n)}\leq C\|f\|_{L^2(S)}?$$
 	 \end{prob}
  It is well known from the  literature that there are only  two types of surfaces namely,  smooth compact surfaces with non-zero Gauss curvature  and  quadratic surfaces, for which the restriction problem, i.e, Problem \ref{problem1}  (or \ref{problem2}), has been settled completely. For smooth compact surfaces with non-zero Gauss curvature,   Stein-Tomas theorem   asserts that   the restriction conjecture  (\ref{Restriction1}) holds for all $1\leq p\leq \frac{2(n+1)}{n+3},$ see \cite{ES, PT}.   On the other hand, for  quadratic surfaces such as paraboloid-like, cone-like, or sphere-like,  Strichartz gave a complete characterization depending on surfaces in \cite{RS}. For a detailed study on restriction problems, we refer to the excellent review  by Tao \cite{tao}.

   The main aim of this article is to validate the Problem \ref{problem1}  or Problem \ref{problem2}  for   the Fourier-Dunkl transform   for certain surfaces and extend it for system of orthonormal functions.  Note that Dunkl transform was introduced by Charles Dunkl (1988) to build a framework for a theory of special functions and integral transforms in several variables related to reflection groups. The basic properties of the Dunkl transforms have been studied by several authors, see \cite{J,CFD,CFD2,T} and the references therein.  Very recently, many authors have been investigating the behavior of the Dunkl transform  to several problems already studied for the Fourier
transform; for instance, Babenko inequality \cite{Bouz1},  uncertainty principles ~\cite{GT,Yakarim}, real Paley-Wiener theorems~\cite{HT}, heat equation ~\cite{MR}, Dunkl Gabor transform \cite{MS}, Dunkl wave equation \cite{H}, Dunkl Schr\"{o}dinger equation \cite{H9,H13},    Dunkl wavelet transform \cite{T3}, Riesz transform \cite{ST3} and so on.

% Note that Dunkl operators were introduced by Charles Dunkl (1989) to build a framework for a theory of special functions and integral transforms in several variables related to reflection groups. Such operators are relevant in physics, namely for the analysis of quantum many-body systems of Calogero-Moser-Sutherland type (see \cite{DV,LV}). From the mathematical analysis point of view, the importance of Dunkl operators lies on the fact that they generalize the theory of symmetric spaces of Euclidean type. There are many developments in harmonic analysis of the operators, which are defined in terms of Dunkl operators in different directions in recent years.  We refer to  \cite{BRS,CFD,CFD2,DX,H,MR}  for detailed study on  Dunkl operators.
   
%------------------------------------------------------------------------------------------------------------------------------------------------------
\subsection{Restriction problem for the Fourier-Dunkl transform}
Let G be a Coxeter group associated with the root system $\mathcal{R}\subset\R^d,~\K$ be the multiplicity function on $\mathcal{R}$ such that $\K(\nu)\geq0,$ for all $\nu\in \mathcal{R}$ and $\mathcal{R}^+$ be a positive root system of $\mathcal{R}.$    Let $h^2_{k}$ and $E_\K$ denotes weight function and Dunkl kernel associated with the root system $\mathcal{R}$, respectively. We refer to   Section \ref{pre2} for  detailed  definitions and   unexplained notations.

%
%Let $h^2_{k}(x):=\prod_{\nu\in R_+}|\langle x, \nu\rangle|^{2\kappa(\nu)}$ be the weight function and $T_\xi$ denotes the corresponding Dunkl operators defined by
%$$T_{\xi}f(x)=\partial_{\xi}f(x)+\sum_{\nu\in R_{+}}\kappa(\nu)\langle\nu, \xi\rangle\frac{f(x)-f(\sigma_{\nu}x)}{\langle\nu, x\rangle},$$
%where $\sigma _\nu$ is the reflection in the hyper plane perpendicular to $\nu,$ and $f\in C'(\R^d).$
%The Dunkl-Laplacian is the second order operator defined by
%$$\Delta_{\kappa}=\sum^d_{j=1}T^2_{e_j},$$
%where $\{e_1, e_2, \ldots e_d\}$ is the standard basis for $\R^d$. Let $E_\K$ denotes the Dunkl kernel.

Let $L^p_{\kappa}(\R^n\times\R^d)$, $1\leq p\leq \infty,$ denotes the space of $L^p$-functions $f(x, y)$ with respect to the measure $h_\kappa^2(y)dxdy,~~\text{for}~~(x, y)\in \R^n\times\R^d.$  %\textcolor{red}{Now we are ready to define the Fourier-Dunkl transform.}
Let  $f\in L^1_\K(\R^n\times\R^d). $ Then  the Fourier-Dunkl transform of $f$, $\widehat f$  is     defined as
\begin{equation}\label{1}
	\widehat f(\xi, \zeta) = \frac{1}{c_\K(2\pi)^{n/2}}\int_{\R^n}\int_{\R^d} f(x, y) E_\K(-i\zeta, y) e^{-ix\cdot \xi} h^2_\K(y) dydx,
\end{equation}
for all $(\xi, \zeta)\in \R^n\times\R^d.$ 

Given a surface $S\subset \R^n\times\R^d,~~n+d\geq 2,$  for all Schwartz class functions $f$,  the restriction operator  is defined as $\rc_Sf =  \widehat f|_S$. Then the restriction problem for the Fourier-Dunkl transform  reads as: 
\begin{prob}\label{problem3}
	For  what values of $p,~~1\leq p < 2,$ do we have
	\begin{equation}\label{5}
		\left(\int_S| \widehat f(\xi, \zeta) |^2 h^2_\K(\xi)d\sigma(\xi,\zeta)\right)^\frac{1}{2}\leq C \|f\|_{L^p_\K(\R^n\times\R^d)}?
	\end{equation}
where $d\sigma(x, y)$ denotes surface measure on $S$ or $(n+d-1)$-dimensional Lebesgue measure endowed with the surface $S.$
\end{prob}
Then the       dual operator of  $\rc_S$,     the extension  operator  $\ecs$     is defined as
\begin{equation}\label{6}
	\ecs f(x, y) = \frac{1}{c_\K (2\pi)^{n/2}}\int_S f(\xi, \zeta) E_\K(i y, \zeta) e^{ix\cdot \xi} h^2_\K(\zeta) d\sigma(\xi, \zeta).
\end{equation}
Thus \eqref{5} is equivalent to the boundedness of $\ecs$ from $L^2(S, h_\K^2(y)d\sigma(x, y))$ to $L^{p'}_\K(\R^n\times\R^d)$, i.e., to say
\begin{equation}\label{ecb}
	\|\ecs f\|_{L^{p'}_\K(\R^n\times\R^d)}\leq C \|f\|_{L^2(S, h_\K^2(y)d\sigma(x, y))},
\end{equation}where $p'$ denotes the conjugate exponent of $p$. Note that,  when $\K\equiv 0,$ the Fourier-Dunkl transform coincides with the Fourier transform on $\R^{n+d}$ and  in which case,  the restriction problem was considered by Strichartz, Stein and P. A. Tomas (as explained at the beginning of this section).

 \subsection{Orthonormal version of restriction problem for the Fourier-Dunkl transform}
 Recently,   a considerable attention has been devoted by several researchers to extend  functional  inequalities from single function to   a system of orthonormal functions in different frameworks. Thus an interesting and natural question  we want to address in this paper is the  generalization of (\ref{ecb})  in the framework of orthonormal systems.   
 More precisely, let $\left(f_{j}\right)_{j \in J}$ be a (possibly infinite) system of orthonormal functions  in $L^2(S, h_\K^2(y)d\sigma(x, y))$, and let $\left(n_{j}\right)_{j \in J} \subset \mathbb{C}$ be a sequence of coefficients,   then one can ask, for which exponents  $1 \leq p \leq 2$, we have
 \begin{align}\label{Intro1}
 	\left\|\sum_{j \in J} n_{j}\left|\mathcal{E}_{S} f_{j}\right|^{2}\right\|_{L^{\frac{p^{\prime}} {2}}_\K(\R^n\times\R^d)}\leq  C\left(  \sum_{j \in J}\left|n_{j}\right|^\alpha\right)^{\frac{1}{\alpha}},
 \end{align}
 for some $\alpha>1$ and  for some positive constant $C$ (independent of $\left(f_{j}\right)$ and $\left(n_{j}\right)$).

 Note that  the idea of extending functional inequalities from single input  to   systems of orthonormal functions is hardly  a new topic. Such generalization involving the orthonormal system is strongly motivated by the theory of  many-body systems in  quantum mechanics, where  a simple description of   $N$ independent fermions particles in the Euclidean space, $\mathbb{R}^{n}$  can be  described by a collection of  $N$ orthonormal functions in $L^{2}\left(\mathbb{R}^{n}\right)$.  For this reason, functional inequalities involving a large number of orthonormal functions are very useful in the mathematical analysis of many body  quantum systems.   The first initiative   work  of such  generalization goes back to the famous work established by Lieb and Thirring \cite{lieb} in  which they extended   versions of certain Gagliardo--Nirenberg--Sobolev inequalities for a single function   and applications to the stability of matter. Later,  in 2013, the    classical Strichartz inequality for the Schr\"odinger propagator $e^{it\Delta}$ \cite{RS} have been substantially generalized  for  a system of orthonormal functions by Frank-Lewin-Lieb-Seiringer \cite{frank} and Frank-Sabin \cite{FS}. We refer to \cite{ben, lee, JP,  MS, shyam} for recent developments in the direction of orthonormal systems.

Motivated by the recent developments, the authors in \cite{JPSH}  proved restriction theorems  for the Fourier-Dunkl  transform and extended it for the system of orthonormal functions for the cone type surface.   Similar to \cite{JPSH}, in this paper, we further extend and   give sufficient condition for the estimates \eqref{5} (or \eqref{ecb})  and find exponents $p$ such that the orthonormal inequality (\ref{Intro1})  holds for certain surfaces in the  Fourier-Dunkl setting mainly using the   ideas of Strichartz \cite{RS} and Frank-Sabin \cite{FS}.

\subsection{Main results}   In this article, we consider  surfaces namely, paraboloid
\begin{equation}\label{sur1}
S_1=\{(x, y)\in\R^n\times\R^d: x_n= x_1^2 + \cdots + x_{n-1}^2-y_1^2- \cdots -y_d^2\},
\end{equation} 
sphere 
\begin{equation}\label{sur2}
S_2=\{(x, y)\in\R^n\times\R^d: |x|^2+|y|^2=1\},
\end{equation}
and two sheeted hyperboloid
\begin{equation}\label{sur3}
S_3=\{(x, y)\in\R^n\times\R^d: |x|^2-|y|^2=1\}.
\end{equation} 
For the paraboloid surface $S_1$,  we get the following  Fourier-Dunkl restriction theorem:
\begin{theorem}[Fourier-Dunkl restriction theorem for paraboloid]\label{FDRTP}
Let $S_1$ be the paraboloid as given in \eqref{sur1},  $n+d\geq 2$, $n\geq 1$, $N_\K= n+d+2\gk$ and $p= \frac{2(N_\K+1)}{N_\K+3}.$ Then the restriction operator  $\rc_{S_1}=\widehat{f}|_{S_1}$ 
%initially defined on $\Sc(\R^n\times\R^d)$ 
can be extended as a bounded operator from $L^p_\K(\R^n\times\R^d)$ to $L^2(S_1, h^2_\K(y)d\si(x, y))$
and we have $$\left(\int_{S_1}| \widehat f(\xi, \zeta) |^2 h^2_\K(\z)d\sigma(\xi,\zeta)\right)^\frac{1}{2}\leq C\|f\|_{L^p_\K(\R^n\times\R^d)},$$  
for all $f\in \Sc(\R^n\times\R^d)$. 
\end{theorem} 
Moreover, using  a   duality principle, we have the following restriction estimates for orthonormal functions.
\begin{theorem}[Restriction estimates for orthonormal functions-paraboloid]\label{REOp}
	Let $S_1$ be the surface as given in \eqref{sur1}, $N_\K= n+d+2\gk\geq 2$ for any (possibly infinite) orthonormal system $f_j$ in  $L^2(S_1, h^2_\K(y)d\si(x, y))$ for any $(n_j)\subset \C,$ we have 
	$$\left\|\sum_j n_j|\ecso f_j|^2\right\|_{L_\K^\frac{N_\K+1}{N_\K-1}(\R^n\times\R^d)}\leq C \left(\sum_j |n_j|^\frac{N_\K+1}{N_\K}\right)^\frac{N_\K}{N_\K+1}.$$  
\end{theorem}

 Further, for the sphere surface $S_2$,  we have the following  Fourier-Dunkl restriction theorem:
\begin{theorem}[Fourier-Dunkl restriction theorem for sphere]\label{FDRTS2}
Let $S_2$ be the unit sphere in $\R^n\times\R^d$ defined  in \eqref{sur2}, $N_\K= n+d+2\gk\geq 2,$ and $1\leq p\leq \frac{2(N_\K+1)}{(N_\K+3)}$. Then  the restriction operator $\rc_{S_2}=\widehat{f}|_{S_2}$ can be extended as bounded operator from $L^p_\K(\R^n\times\R^d)$ into $L^2(S_2, h^2_\K(y)d\si(x, y))$ and we have 
$$\left(\int_{S_2}| \widehat f(\xi, \zeta) |^2 h^2_\K(\z)d\sigma(\xi,\zeta)\right)^\frac{1}{2}\leq C\|f\|_{L^p_\K(\R^n\times\R^d)},$$     
for all $f\in \Sc(\R^n\times\R^d)$.
\end{theorem}
\begin{rem}
When $n=0$, Fourier-Dunkl transform becomes Dunkl transform on $\mathbb{R}^d$ and Theorem \ref{FDRTS2} gives the restriction theorem for the Dunkl transform on unit sphere for general root system with non-negative multiplicity  function $\K$.  Thus Theorem \ref{FDRTS2} generalizes    Corollary 1.3 of Dai and Ye \cite{dai}, where the authors proved  restriction theorem for the Dunkl transform on unit sphere  associated with the reflection group $\mathbb{Z}_2^d$.
\end{rem}

Using  the same   duality principle, we have the following restriction estimates for orthonormal functions.
\begin{theorem}[Restriction estimates for orthonormal functions-sphere]\label{REOs}
	Let $S_2$ be the unit sphere in $\R^n\times\R^d$ given in \eqref{sur2}, $N_\K= n+d+2\gk\geq 2.$ Then for any (possibly infinite) orthonormal system $(f_j)$ in  $L^2(S_2, h^2_\K(y)d\si(x, y))$ for any $(n_j)\subset \C,$ we have 
	$$\left\|\sum_j n_j|\ecstw f_j|^2\right\|_{L_\K^r(\R^n\times\R^d)}\leq C \left(\sum_j |n_j|^{\frac{2r}{r+1}}\right)^{\frac{r+1}{2r}},$$
	for $\frac{N_\K+1}{N_\K-1}\leq r<\infty$.
\end{theorem}

 Finally, for the hyperboloid surface $S_3$,  we have the following  Fourier-Dunkl restriction theorem:
\begin{theorem}[Fourier-Dunkl restriction theorem for hyperboloid]\label{FDRTh}
Let $S_3$ be the surface in $\R^n\times\R^d$ given in (\ref{sur3}) with $n\neq 0$ and $N_\K= n+d+2\gk\geq 2.$

\textbf{Case $1$.} Let $n=d=1$ and $\gk=0$ so that $N_\K=2.$ Then the Fourier-Dunkl transform becomes Fourier transform on $\R^2$ and for  $1< p\leq \frac{6}{5},$  the restriction operator $\rc_{S_3}=\widehat{f}|_{S_3}$ can be extended as bounded operator from $L^p_0(\R\times\R)$ into $L^2(S_3, h^2_\K(y)d\si(x, y))$ and we have 
$$\left(\int_{S_3}| \widehat f(\xi, \zeta) |^2 d\sigma(\xi,\zeta)\right)^\frac{1}{2}\leq C\|f\|_{L^p_0(\R\times\R)},$$ 
for all $f\in \Sc(\R\times\R).$

\textbf{Case $2$.}  Let $n\geq1,~~ d\geq1, N_\K>2$, and $\frac{2N_\K}{N_\K+2}\leq p\leq \frac{2(N_\K+1)}{(N_\K+3)}.$ Then  the restriction operator $\rc_{S_3}=\widehat{f}|_{S_3}$ can be extended as bounded operator from $L^p_\K(\R^n\times\R^d)$ into $L^2(S_3, h^2_\K(y)d\si(x, y))$    and we have
$$\left(\int_{S_3}| \widehat f(\xi, \zeta) |^2 h^2_\K(\z)d\sigma(\xi,\zeta)\right)^\frac{1}{2}\leq C\|f\|_{L^p_\K(\R^n\times\R^d)},$$ 
for all $f\in \Sc(\R^n\times\R^d)$.
\end{theorem}
  Again, using  the same  duality principle, we have the following restriction estimates for orthonormal functions.
\begin{theorem}[Restriction estimates for orthonormal functions-hyperboloid]\label{REOh}
	Let $S_3$ be the unit sphere in $\R^n\times\R^d$ given in \eqref{sur3}, $N_\K= n+d+2\gk\geq 2.$ Then 
	
	Case $1.$ If $n= d= 1$ and $\gk=0$ so that $N_\K=2,$ then for any (possibly infinite) orthonormal system $(f_j)$ in $L^2(S_3, d\si(x, y))$ and for any $(n_j)\subset \C,$ we have    
	$$\left\|\sum_j n_j|\ecsth f_j|^2\right\|_{L_\K^r(\R\times\R)}\leq C \left(\sum_j |n_j|^{\frac{2r}{r+1}}\right)^{\frac{r+1}{2r}},$$
	for $3\leq r<\infty$.
	
	Case $2.$ If $n\geq 1,~~d\geq 1$ and $N_\K>2$, then for any (possibly infinite) orthonormal system $(f_j)$ in $L^2(S_3, h^2_\K(y)d\si(x, y))$ and for any $(n_j)\subset \C,$ we have    
	$$\left\|\sum_j n_j|\ecsth f_j|^2\right\|_{L_\K^r(\R^n\times\R^d)}\leq C \left(\sum_j |n_j|^{\frac{2r}{r+1}}\right)^{\frac{r+1}{2r}},$$
	for $\frac{N_\K+1}{N_\K-1}\leq r \leq\frac{N_\K}{N_\K-2}$.
\end{theorem}

As we earlier discussed,   Strichartz inequalities are important applications of the restriction estimates of quadratic surfaces,  which are useful tools to study non-linear Schr\"odinger equations.   For instance, if we consider the paraboloid surface  $S=\{(\om, \z)\in\R\times\R^d: \om= -|\z|^2\},$ then  for a particular choice of $f\in L^1(S, d\mu),$  where $d\mu(\om, \z)=d\z,$ we get 
\begin{align*}
	\ecs f(t, y)=\frac{1}{\sqrt{2\pi}}(e^{it\Dk}\phi)(y),  
\end{align*}     
for    some $\phi: \R^d\to \C,$ such that $f(\om, \z)= \mathcal{F}_\K\phi(\z)$. In this case we get the following  orthonormal Strichartz inequalities for the Schr\"odinger  propagator $e^{it \Dk}$ associated with Dunkl Laplacian operator:
\begin{theorem}[Strichartz estimates for orthonormal functions for Schr\"odinger  propagator] 
	Let $d\geq 1$ and $p, q\geq 1$ such that $\frac{2}{p}+\frac{d+2\gk}{q}=d+2\gk,~~ 1\leq q<1+\frac{2}{d+2\gk-1}.$ Then for any (possible infinite) orthonormal system $\phi_j$ in $L^2(\R^d, h^2_\K(y) dy)$ and for any $(n_j)\subset \C$ we have 
	\begin{align}\label{1.1}
		\Big\|\sum_j n_j|e^{it\Dk}\phi_j|^2\Big\|_{L^p_t(\R, L^q_y(\R^d, h^2_\K(y) dy))}\leq C\Big(\sum_j |n_j|^\frac{2q}{q+1}\Big)^\frac{q+1}{2q},
	\end{align} 
with $C>0$ independent of $(n_j)$ and $(\phi_j).$
\end{theorem}
\begin{rem}
	 Here we want to make a note that the inequality (\ref{1.1}) is proved  in \cite{JP, MS} by means of similar estimates for Dunkl-Hermite operators and Schr\"odinger kernel relation, see Theorem 6.1 of \cite{MS} or  Lemma 3.1 of \cite{JP}.

	   %using the classical Strichartz estimates for the free Schr\"odinger propagator $e^{-i t \Delta_{k, 2}} $  for  orthonormal systems of initial data and the kernel relation  between the semigroups $e^{-i t \Delta_{k, a}}$ and $e^{i \frac{t}{2} \Delta_{k}},$ we prove   Strichartz estimates for orthonormal systems of initial data    associated with  the Dunkl operator $ \Delta_k $ on $\mathbb{R}^n$.  Finally,  we present some applications  to our   aforementioned results.
	 
\end{rem}

%	Further, 

Again,  if we  consider the surface  $S=\{(\om, \z)\in\R\times\R^d: \om^2=1+|\z|^2\}$ with the measure $d\mu(\om, \z)=\frac{d\z}{2\sqrt{1+|\z|^2}},$ then    if we choose $f(\om, \z)= 2\hspace{0.2cm}\mathds{1}_{(\om>0)} \sqrt{1+|\z|^2} \mathcal{F}_\K(\phi)(\z),$ then we obtain
$$\ecs f(t, y)=\frac{1}{\sqrt{2\pi}}e^{it\sqrt{1-\Dk}}\phi(y),\hspace{0.5cm}\forall (t, y)\in(\R\times\R^d).$$
 In this case, we get the corresponding  orthonormal Strichartz inequalities for the Klein-Gordon  propagator $e^{it \sqrt{1-\Dk}}$ associated with Dunkl Laplacian operator in the following result:
\begin{theorem}[Strichartz estimates for orthonormal functions for the Klein-Gordon  propagator] 
	Assume $d\geq 1.$ Let $1+\frac{2}{d+2\gk}\leq r\leq 1+\frac{2}{d+2\gk-1}$ if $d+2\gk>1$ and $3\leq r<\infty$ if $d+2\gk=1.$ For any (possibly infinite) orthonormal system $(\phi_j)$ in $H_\K^{1/2}(\R^d)$ and for any $(n_j)\subset \C,$ we have 
	$$\Big\|\sum_j n_j|e^{it\sqrt{1-\Dk}}\phi_j|^2\Big\|_{L^r(\R\times\R^d)}\leq C \Big(\sum_j|n_j|^\frac{2r}{r+1}\Big)^\frac{r+1}{2r},$$
	with $C>0$ independent of $(n_j)$ and $(\phi_j).$
\end{theorem}
 An immediate consequences of the above result is as follows. 
\begin{cor}
	Let $d\geq 1,$ suppose $1+\frac{2}{d+2\gk}\leq r_0\leq 1+\frac{2}{d+2\gk-1},$ for $1\leq r\leq r_0,~~q\geq 1$ and $s\geq 0$ such that $\frac{1}{q}=\frac{1}{1-r_0}\big(1-\frac{1}{r}\big)$ and $s=\frac{r_0}{r_0-1}\big(\frac{1}{2}-\frac{1}{2r}\big).$ Then for all families of orthonormal functions $(\phi_j)$ in $H^s_\K(\R^d),$ we have that 
	\begin{equation}\label{kgec}
		\Big\|\sum_j n_j|e^{it\sqrt{1-\Dk}}\phi_j|^2\Big\|_{L^q_t(\R, L^r_y(\R^d, h^2_\K(y) dy)}\leq C \Big(\sum_j|n_j|^\bt\Big)^\frac{1}{\bt},
	\end{equation}
	where $\bt= \frac{2r}{r+1}.$
\end{cor}

Apart from introduction, the organization of the article is as follows: 
\begin{itemize}
	\item In Section \ref{pre2},  we   recall harmonic analysis related to Dunkl operators, Fourier-Dunkl transform and their certain properties.
	\item In Section \ref{gtrt},  we  recall crucial ideas to prove restriction theorem for a general surface and its  extension to  an orthonormal family of functions.
	
	\item In   Section \ref{sec3}, we   find the Fourier-Dunkl transform of certain generalized functions which are required to prove the main result. Particularly,  the Fourier-Dunkl transform of   $(x_n-|x'|^2+|y|^2)_+^z,~~(1-|x|^2-|y|^2)_+^z$, and $(1-|x|^2+|y|^2)_+^z,$ where $x'= (x_1, \ldots, x_{n-1}).$
	
	\item In   Section \ref{pmr}, we prove our  main results of this paper.
	
	\item In   Section \ref{sec5}, 	we  establish   Strichartz estimates for orthonormal families of initial data associated  with  Dunkl Laplacian propagator $e^{it\Dk}$  and Klein-Gordon propagator $e^{it\sqrt{1-\Dk}}$. 
\end{itemize}

%In Section \ref{pre2}, we recall the general Dunkl setting, Schatten spaces and more about the  Fourier-Dunkl transform. In Section \ref{gtrt}, we develop the similar method as in \cite{RS} to prove the restriction theorem then give a sufficient condition in terms of generalized functions and we deduce the Strichartz inequality for the orthonormal functions associated with square root of Dunkl Laplacian as an application of restriction theorem for the cone-hyper-surface. We compute the Fourier-Dunkl transform of certain generalized functions which are useful in proving the main results in Section \ref{fdtcgf}. In Section \ref{pmr} we give the proofs of our main results.   Finally, we shall write $C$ to denote positive constants independent of significant quantities the meaning of which can change from one occurrence to another.

%------------------------------------------------------------------------------------------------------------------------------------------------------

%-----------------------------------------------------------------------------------------------------------------------------------------------------------------------

\section{Preliminaries}\label{pre2}
In this section, we recall harmonic analysis related to Dunkl operators, Fourier-Dunkl transform and their certain properties   which we are going to use    throughout this article.    A complete account of  harmonic analysis  related to  Dunkl operators  can be found in \cite{BRS,CFD,CFD2,DX,H,MR,ST2,T}. However, we mainly adopt the notation and terminology given in \cite{JPSH}.

\subsection{General Dunkl setting}
 The basic ingredients in the theory of Dunkl operators are the root systems and finite reflection groups associated with them. For $\nu\in\R^d\setminus\{0\}$, we denote by $\sigma _\nu$, the reflection in the hyperplane perpendicular to $\nu,$ i.e.,
$$\sigma _\nu(y)=y-2\frac{\langle\nu, y\rangle}{|\nu|^2}\nu.$$

Let $O(d)$ be the group of orthogonal matrices acting on $\R^d.$ Given a root system $\mathcal{R}$, we associate a finite subgroup $G\subset O(d)$, the reflection group which is generated by the reflections $\{\sigma_{\nu}: \nu\in \Rc \}.$ A function $\kappa : \Rc\rightarrow \C$ is said to be a multiplicity function on $\Rc,$ if it is invariant under the natural action of $G$ on $\Rc,$ i.e., $\kappa(g\nu)=\kappa(\nu) $ for all $\nu\in \Rc$ and $g\in G.$

Every root system can be written as a disjoint union $\Rc = \Rc_+ \cup (-\Rc_+)$, where $\Rc_+$ and $-\Rc_+$ are separated by a hyperplane through the origin. Such $\Rc_+$ is the set of all positive roots in $\Rc$. Of course, its choice is not unique.

The weight function associated with the root system $\mathcal{R}$ and the multiplicity function $\kappa$ is defined by
$$h^2_{k}(y):=\prod_{\nu\in \Rc_+}|\langle y, \nu\rangle|^{2\kappa(\nu)}.$$
Note that $h^2_{\kappa}(y)$ is $G$-invariant and homogenous of degree $2\gamma_\kappa$, where by definition
$$\gamma_\kappa:=\sum_{\nu\in \Rc_{+}}\kappa(\nu).$$
% Let $L^p_{\kappa}(\R^d)$, $1\leq p\leq \infty,$  stands for the space of $L^p$- functions with respect to the measure $h_\kappa^2(x) dx$, and $L^p_t L^q_{\kappa, x}(\T \times \R^d), 1\leq p, q \leq \infty$ stands for the space of all measurable functions $h(t, x)$ on $\T  \times \R^d$   for which
%$$\|h\|_{L^p_t L^q_{\kappa, x}(\T \times \R^d)}:=\|\|h(t, \cdot)\|_{L^q_{\kappa}(\R^d)}\|_{L^p(\T)}<\infty.$$
%We may consider $\T=(-\pi, \pi), \R$ or any interval in $\R$ with Lebesgue measure. $L^p_\K(\T\times \R^d):=L^p_t L^p_{\kappa, x}(\T \times \R^d)$.

Now we recall the difference-differential operators, introduced and studied by C. F. Dunkl (for $\kappa\geq 0$) see \cite{CFD, CFD2}. These operators are also called Dunkl operators,  which are also the analogues (generalizations) of directional derivatives. We fix a root system $\mathcal{R}$ with a positive subsystem $\Rc_+$ and the associated reflection group $G.$ We also fix a nonnegative multiplicity function $\kappa$ defined on $\Rc.$

For $\z \in\R^d,$ the Dunkl-operator $T_{\z}:=T_{\z}(\kappa)$ is defined by
$$T_{\z}f(y)=\partial_{\z}f(y)+\sum_{\nu\in \Rc_{+}}\kappa(\nu)\langle\nu, \z\rangle\frac{f(y)-f(\sigma_{\nu}y)}{\langle\nu, y\rangle}$$ for smooth functions $f$ on $\R^d.$ Here $\partial_{\z}$ denotes the directional derivative along $\z.$ For the standard coordinate vectors $\z=e_{j}$ of $\R^d$ we use the abbreviation $T_{j}=T_{e_j}.$

Let $\mathscr{P}$ be the space of all polynomials with complex coefficients in $d$-variables and $\mathscr{P}_m$ be the subspace of homogeneous polynomials of degree $m$. The Dunkl-operators $T_\z$ and directional derivatives $\partial_\z$ are closely related and intertwined by an isomorphism on $\mathscr{P}.$ Indeed, if the multiplicity function $\kappa$ is non-negative then by Theorem 2.3 and Proposition 2.3 in R\"osler \cite{MR}, there exist a unique linear isomorphism (intertwining operator) $V_\kappa$ of $\mathscr{P}$ such that $V_\kappa(\mathscr{P}_m)=\mathscr{P}_m, V_\kappa|_{\mathscr{P}_0}=id$ and $T_{\z} V_\kappa=V_\kappa\partial_\z$ for all $\z\in\R^d.$ It can be checked that $V_\kappa\circ g=g\circ V_\kappa$ for all $g\in G.$ For $\z\in\C^d,$ define
$$E_\kappa(y, \z):=V_\kappa\left(e^{\langle \cdot, \z\rangle}\right)(y),~ y\in\R^d.$$
The function $E_\kappa$ is called the Dunkl-kernel, or $\kappa$-exponential kernel, associated with $G$ and $\kappa,$ see \cite{CFD4}.  We collect some important properties of this kernel in the following proposition.
\begin{prop} Let $y, \z \in \C^d$, $\la \in \C $ and $g \in G $,
	\begin{enumerate}
		\item $E_{\K}(y, \z)=E_{\K}(\z, y)$.
		\item $E_{\K}(\la y, \z)=E_{\K}(y, \la \z)$ and $E_{\K}(g y, g \z)=E_{\K}(y, \z)$.
		\item $\overline{E_{\K}(y, \z)}=E_{\K}(\bar{y},  \bar{\z})$.
		\item $|E_{\K}( y, \z)|\leq e^{|y|\:|\z|}$; \; \;  $|E_{\K}( iy, \z)|\leq 1 $ if $y, \z \in \R^d.$
	\end{enumerate}
\end{prop}

For proof of the above proposition, we refer to R\"{o}sler \cite{MR}. Now we define the Dunkl transform, a generalization of the Fourier transform. For $1\leq p <\infty$, let
$L^p(\R^d, h_{\K}^2dy)$ be the space of all measurable functions $f$ defined on $\R^d$ for which the norms
$$\|f \|_{\K,p}=\bigg(\int_{\R^d}|f(y)|^p h_{\K}^2(y)\,dy \bigg)^{\tfrac{1}{p}} $$ are finite.
The space $L^{\infty}(\R^d, h_{\K}^2 dy)$ is defined
in the usual way. The Dunkl transform $\mathcal{F}_\K f$ of $f\in L^1(\R^d,
h_{\K}^2dy)$ is given by

\be \label{eq:Dunkltr} \mathcal{F}_\K
f(\z) = \frac{1}{c_\K}\int_{\R^d} f(y)E_\K(-i\z, y) h_\K^2(y) dy,\ee where
\be \label{ck}
c_\K=\int_{\R^d} e^{-|y|^2/2} h_\K^2(y) dy.
\ee
The Dunkl transform shares many important properties with the Fourier transform. For example, we have the Plancherel theorem \be \label{eq:Dunklplanch} \int_{\R^d} |\mathcal{F}_\K f(\z) |^2 h_\K^2(\z) d\z =  \int_{\R^d} |f(y)|^2  h_\K^2(y) dy \ee for
all $ f \in L^1 \cap L^2(\R^d, h_\K(y)^2 dy) $ and the inversion formula \be \label{eq:Dunklinv} f(y) = \frac{1}{c_\K}  \int_{\R^d}
\mathcal{F}_\K f(\z) E_\K(i y,\z)  h_\K^2(\z) d\z \ee for all $ f \in L^1(\R^d, h_\K^2 dy) $ provided $ \mathcal{F}_\K f $ is also
in $L^1(\R^d, h_\K^2(y) dy).$ The Dunkl transform was introduced in \cite{CFD5} for non-negative multiplicity functions and further
studied in \cite{DJ} in the more general case $\Re{(\K)}\geq 0.$ We refer to these two papers for more details about the Dunkl transform.

For given $y\in \R^d$, the Dunkl translation (generalized translation) $\tau_y$ is an operator on $L^2(\R^d, h_{\K}^2dy)$ defined by the equation
\be \mathcal{F}_\K (\tau_y f)(\z)=E_\K(-i\z, y)\mathcal{F}_\K f(\z),~~~\forall \z, y\in\R^d.\ee
From the definition it can be seen that $\tau_y f(\z)=\tau_{-\z} f(-y)$. It is still an open problem whether the Dunkl translation can be extended to a bounded operator on $L^p(\R^d, h_\K^2(y) dy)$ for any $p\geq 1$. But the Dunkl translation can be extended to $L^p_{rad}(\R^d, h_\K^2(y)dy)$ for $1\leq p \leq \infty$ as a bounded operator and it can be seen by using the expression $$\tau_y(f)(\z)=\int_{\R^n} f_0\big(\sqrt{|\z|^2+|y|^2+2\La\z, w\Ra}\big)d\nu_y(w),$$ where $f(y)=f_0(|y|)$ and $d\nu_y$ is a  probability measure supported in convex hull of $\{gy~:~g\in G\}$, see \cite{ST2,MR}.

Let $\Sc(\R^d)$ denote the space of all Schwartz class functions on $\mathbb{R}^d$. The following important property of the Dunkl translation is very useful.
\be \label{transp} \int_{\R^d} \tau_y f(\z) g(\z) h_\K^2(\z) d\z=\int_{\R^d}f(\z)  \tau_{-y} g(\z) h_\K^2(\z) d\z \ee for all $f, g\in \Sc(\R^d)$. Using the Dunkl translation, we define the Dunkl convolution of functions $f, g\in \Sc(\R^d)$ by
\begin{equation}\label{dc}
	f\ast_\K g(y)=\int_{\R^d}\tau_y \tilde{f}(\z) g(\z) h^2_\K(\z) d\z,
\end{equation} where $\tilde f(y)=f(-y)$.
We note the following important properties of Dunkl convolution
\be \mathcal{F}_\K (f\ast_\K g)=\mathcal{F}_\K (f)  \mathcal{F}_\K( g)\ee for all $f, g\in \Sc(\R^d)$
and
\be \|f \ast_\K g\|_{\K, \infty}\leq \|g\|_{\K, \infty} \|f\|_{\K, 1} \ee
for all radial functions $g$ and $f\in\Sc(\R^d)$. We refer to the  papers \cite{ST2,T} for more details about the Dunkl covolution.

The Dunkl Laplacian is the second-order operator defined by
$$\Delta_{\kappa}=\sum^d_{j=1}T^2_j$$ which can be explicitly calculated; see Theorem 4.4.9 in Dunkl-Xu \cite{DX}. It can be seen that $\Delta_\kappa=\sum^d_{j=1}T^2_{\xi_j}$ for any orthonormal basis $\{\xi_1, \xi_2,\ldots, \xi_d\}$ of $\R^d$, see \cite{CFD2}. In \cite{AH}, the authors proved the operator $-\Dk$ is essentially a self-adjoint positive operator on $L^2(\R^d, h_\K^2(y)dy)$ and the Dunkl transform gives its spectral decomposition. For any measurable function $m$ on $\R$, the operator $m(-\D_k)$ defined by
$$m(-\D_k) f(y)=\frac{1}{c_\K}\int_{\R^d}m(|\z|^2) \mathcal{F}_\K(f)(\z)  E_\K(iy, \z)h_\K^2(\z) d\z.$$ It is easy to see that $m(-\D_k)$ is a bounded operator on $L^2(\R^d, h_\K^2(y)dy)$ if and only if $m$ is bounded. Moreover,  $e^{it\sqrt{1-\Dk}}$ and $e^{it\Dk}$ are unitary operators for any $t \in \R$ which are defined by
$$e^{it\sqrt{1-\Dk}}f(y)=\frac{1}{c_\K}\int_{\R^d} e^{it\sqrt{1+|\z|^2}} \mathcal{F}_\K(f)(\z)  E_\K(iy, \z)h_\K^2(\z) d\z.$$
and
$$e^{it\Dk}f(y)=\frac{1}{c_\K}\int_{\R^d} e^{-it|\z|^2} \mathcal{F}_\K(f)(\z)  E_\K(iy, \z)h_\K^2(\z) d\z.$$
For $s\geq 0,$ we denote  $ {H}^s_\K(\R^d)$ for  the completion of $\Sc(\R^d)$ with respect to the norm given by $$\|\phi\|^2_{{H}^s_\K(\R^d)}=\int_{\R^d} (1+|\z|^{2})^s |\mathcal{F}_\K\phi(\z)|^2 h^2_\K(\z)d\z.$$ 
The space  $ {H}^s_\K(\R^d)$  is known as Dunkl-Sobolev space. For a detailed study about the space $ {H}^s_\K(\R^d)$, we refer to \cite{M2}.

%------------------------------------------------------------------------------------------------------------------------------------------------------

\subsection{Schatten spaces}\label{Sbss}
Let $\Hc$ be complex separable Hilbert space and $A$ be a compact operator on $\Hc$.  We say that $A\in \s^p(\Hc),$ Schatten space for $ 1\leqslant p<\infty,$ if $Tr|A|^p<\infty,$ where $|A|=\sqrt{A^*A}$ and $Tr A$ is trace of $A$. For $A\in\s^p(\Hc),$ the Schatten $p$-norm of $A$ in $\s^p(\Hc)$ is defined by
$$\|A\|_{\s^p(\Hc)}=(Tr|A|^p)^\frac{1}{p}.$$
It can be verified that  $\|A\|_{\s^p(\Hc)}=\left(\sum_{j}|\lambda_{j}|^p\right)^\frac{1}{p},$ where $\lambda_1\geqslant\lambda_2\geqslant \dots \geqslant\lambda_n\geqslant \dots \geqslant 0$ are the singular values of $A,$ i.e.,  the eigenvalues of $|A|=\sqrt{A^*A}.$

When $p=2,$  the Schatten $p$-norm coincides with the Hilbert-Schmidt norm.  Also when $p=\infty,$ we define $\|A\|_{\s^\infty(\Hc)}$ to be the operator norm of $A$ on $\Hc.$ For more details of Schatten spaces and its applications on some problems of harmonic analysis, we refer  to \cite{S,W}.

\subsection{Revisit of Fourier-Dunkl transform}
Let $\Sc(\R^n \times \R^d)$ denotes the space of all Schwartz class functions on $\R^n\times \R^d$ and for  $f\in\Sc(\R^n\times\R^d)$, $\widehat{f}^0(\xi, y)$ denotes the Fourier transform of $f$ in first variable, i.e.,
\begin{equation}\label{ft}
	\widehat{f}^0(\xi, y):=\frac{1}{(2\pi)^{n/2}}\int_{\R^n}f(x, y)e^{-ix\cdot\xi}dx,
\end{equation}
and ${\widehat{f}}^\K(x, \z)$ denotes the  Dunkl transform of $f$ in second variable, i.e.,
\begin{equation}\label{dt}
	{\widehat{f}}^\K(x, \z):=\frac{1}{c_\K}\int_{\R^d}f(x, y)E_\K(-i\z, y)h^2_\K(y)dy,
\end{equation}
where $c_\K$ is defined as in $\eqref{ck}$.
%Let $f\in L^1_\K(\R^n\times\R^d),~~~\widehat{f}$ denotes Fourier-Dunkl transform and is defined by
Let  $f\in\Sc(\R^n\times\R^d).$ Then recall that the Fourier-Dunkl transform of $f$ is denoted by $\widehat{f}$ and is defined by
$$\widehat{f}(\xi, \zeta) := \frac{1}{c_\K(2\pi)^{n/2}}\int_{\R^n}\int_{\R^d} f(x, y) E_\K(-i\zeta, y) e^{-ix\cdot \xi} h^2_\K(y) dydx,~~~(\xi, \z)\in \R^n\times \R^d.$$ The inverse Fourier-Dunkl transform is defined by $\widecheck{f}(\xi, \z)=\widehat{f}(-\xi, -\z).$
%by using \eqref{ft} and \eqref{dt} we will define the Fourier-Dunkl transform and denote by $\widehat{f}$ and is defined by
%\begin{equation}\label{fdt}
%\widehat{f}(\xi, \zeta) := \frac{1}{(2\pi)^{n/2}}\int_{\R^n}\int_{\R^d} f(x, y) E_\K(-i\zeta, y) e^{-ix\cdot \xi} h^2_\K(y) dydx,
%\end{equation}
Thus by \eqref{ft} and \eqref{dt}, we can write that
\begin{equation}
	\widehat{f}(\xi, \z)=\int_{\R^n}\widehat{f}^{\K}(x, \z)e^{-ix\cdot\xi}dx= (\widehat{f}^{\K})^{\wedge_0}(\xi, \z)=(\widehat{f}^0)^{\wedge_\K}(\xi, \z).
\end{equation}
\begin{rem}
	For a given root system $\mathcal{R}$ on $\R^d$ and a multiplicity function $\K$ on it, let  $\widetilde{\mathcal{R}}:=\{(0, \nu)\in \R^n\times \R^d: \nu \in \mathcal{R}\}$ and $\widetilde{\K}(0, \nu):=\K(\nu)$ for all $\nu \in  \mathcal{R}$, then      $\widetilde{\mathcal{R}}$ is  also a root system on $\R^{n+d}$ and $\widetilde{\K}$ is a multiplicity function on $\widetilde{ \mathcal{R}}$. As a result, the Fourier-Dunkl transform associated with the root system $\mathcal{R}$ and multiplicity function $\K$ is nothing but Dunkl transform on $\R^{n+d}$ associated with the root system  $\widetilde{\mathcal{R}}$ and the multiplicity function $\tilde{\K}$.
\end{rem}
%\begin{theorem}[Plancherel theorem]
%Let $f\in L^2_\K(\R^n\times\R^d).$ Then $$\|f\|^2_{ L^2_\K(\R^n\times\R^d)}=\|\widehat f\|^2_{ L^2_\K(\R^n\times\R^d)}.$$
%\end{theorem}
%\begin{proof}
%
%\end{proof} We have the following essential properties of the Dunkl-Fourier transform that can be proven easily.
%\begin{lem}
	Let $f\in\Sc(\R^n\times\R^d).$ Then for $j=1, 2,\ldots, n$ and $k=1, 2,\ldots, d$, we have the the  following essential properties  of Fourier-Dunkl transform:
	\begin{enumerate}
		\item $(x_jf\widehat{)}(\xi, \z)=i\frac{\partial}{\partial \xi_j}\widehat f(\xi, \z),$
		\item $(y_kf\widehat{)}(\xi, \z)=iT_k \widehat f(\xi, \z),$
		\item $(\frac{\partial f}{\partial x_j}\widehat{)}(\xi, \z)=i\xi_j \widehat f(\xi, \z),$
		\item $(T_kf\widehat{)}(\xi, \z)=i\z_k  \widehat f(\xi, \z),$
		\item $\widehat{f}\in \Sc(\R^n \times \R^d)$.
	\end{enumerate}
%\end{lem}
It can be verified that
\begin{equation}\label{2}
	\|\widehat f\|_{L^2_\K(\R^n\times\R^d)} = \|f\|_{L^2_\K(\R^n\times\R^d)}
\end{equation}
and
\begin{equation}\label{3}
	f(x, y) = \frac{1}{c_\K (2\pi)^{n/2}}\int_{\R^n}\int_{\R^d} \widehat f(\xi, \zeta) E_\K(i y, \zeta) e^{ix\cdot \xi} h^2_\K(\zeta) d\zeta d\xi, \quad (x, y)\in \R^n\times\R^d.
\end{equation}
Inverse Fourier-Dunkl transform is given by
\begin{equation}\label{4}
	\widecheck{F}(\xi, \zeta) = \widehat F(-\xi, -\zeta).
\end{equation}
Let $f, g\in\Sc(\R^n\times\R^d),$ then  $f\ast g$ denotes the convolution of $f$ and $g$ in first variable, i.e.,
\begin{equation}\label{fc}
	f\ast g(x, y):=\int_{\R^n}f(x-\xi, y) g(\xi, y)d\xi,
\end{equation}
and $f\ast_Dg$ denotes  the Dunkl convolution of $f$ and $g$ in second variable, i.e.,
\begin{equation}\label{dc2}
	f\ast_Dg(x, y)=\int_{\R^d}\tau_y \tilde{f} (x, \z) g(x, \z) h^2_\K(\z) d\z,
\end{equation}
where $\tau_y$ is  the Dunkl translation (generalized translation) in second variable and $\tilde{f}(x, y)=f(x, -y)$. 
%\begin{defn}

	Again for $f, g\in\Sc(\R^n\times\R^d).$ The Fourier-Dunkl convolution of $f$ and $g$ is denoted by $f\ast_{FD}g$ and is defined by
	\begin{align}\label{FD}
		f\ast_{FD}g(x, y)=\int_{\R^n}\int_{\R^d}\tau_y\tilde{f}(x-\xi, \z)g(\xi, \z)h^2_\K(\z)d\z d\xi.
	\end{align}
	By  \eqref{fc} and \eqref{dc}, we also can write   the Fourier-Dunkl convolution as
	\begin{align}
		f\ast_{FD}g(x, y)&=\int_{\R^n}\big(f(x-\xi, \cdot)\ast_D g(\xi, \cdot)\big) (y)d\xi,\nonumber\\
		f\ast_{FD}g(x, y)&=\int_{\R^d}\big(\tau_y\tilde{f}(\cdot, \z)\ast g(\cdot, \z)\big)(x)h^2_\K(\z)d\z.
	\end{align}
Moreover, we have 
\begin{align}\label{lem:con}
 (f\ast_{FD} g\widehat{)}(\xi, \z)=\widehat{f}(\xi, \z)\widehat{g}(\xi, \z),
\end{align}
and  if $g$ is radial in second variable, then
%\begin{proof}From the definition (\ref{1}),  we have
%	\begin{align*}
%		(f\ast_{FD} g\widehat{)}(\xi, \z)&=\int_{\R^n}\int_{\R^d}{(f\ast_{FD} g)}(x, y)E_\K(-i\z, y)e^{-ix\cdot \xi}h^2_{\K}(y)dy dx\\
%		&=\int_{\R^d}\int_{\R^d}(\tau_y\tilde{f}(\cdot, \z)\ast g(\cdot, \z)\widehat{)^{0}}(\xi)h^2_\K(\z)d\z E_\K(-i\z, y)h^2_\K(y)dy\\
%		&=\int_{\R^d}\int_{\R^d}\tau_y\widehat{\tilde{f}}^{0}(\xi, \z) \widehat{g^{0}}(\xi, -\z) h^2_\K(\z) d\z E_\K(-i\z, y) h^2_\K(y) dy\\
%		&=\int_{\R^d}\widehat{f}^{0}\ast_D \widehat{g^{0}}(\xi, y)E_\K(-i\z, y)h^2_\K(y)dy\\
%		&=(\widehat{f}^{0}\ast_D \widehat{g^{0}}\widehat{) }^{\K}(\xi, \z)\\
%		&=(\widehat{f}^0)^{\wedge_\K}(\xi, \z) (\widehat{g^{0}})^{\wedge_\K}(\xi, \z)\\
%		&=\widehat{f}(\xi, \z)\widehat{g}(\xi, \z).
%	\end{align*}
%	This completes the proof of the lemma.
%\end{proof}
%\begin{rem}
%If $T\in\Sc'(\R^n\times\R^d)$ and $g$ be a Schwartz class function on $\R^n\times\R^d$ then $$(T\ast_{FD}g\widehat{)}=\widehat{T}\widehat{g}.$$
%\end{rem}
\begin{align}\label{prop:young}
  \|f \ast_{FD} g\|_{L^\infty_\K(\R^n\times\R^d)}\leq\|f\|_{L^1_\K(\R^n\times\R^d)}\|g\|_{L^\infty_\K(\R^n\times\R^d)}.
\end{align}
%\begin{proof}
%	Using the   fact that $\|\tau_yg\|_{L^\infty_\K(\R^n\times\R^d)}\leq  \|g\|_{L^\infty_\K(\R^n\times\R^d)}$, we have
%	\begin{align*}
%		|(f\ast_{FD}g)(x, y)|&\leq \int_{\R^n}\int_{\R^d}|f(\xi, \z)||\tau_{y}\tilde{g}(x-\xi, \z)| h_\K^2(\z) d\z d\xi\\
%		& \leq \|g\|_{L^\infty_\K(\R^n\times\R^d)}\|f\|_{L^1_\K(\R^n\times\R^d)}
%	\end{align*}
%	for all $(x, y)\in \R^n \times \R^d$.
%\end{proof}
\begin{defn}\cite{JPSH}
	Let $T\in \Sc'(\R^n\times \R^d)$ be a tempered distribution. Then its Fourier-Dunkl transform $\widehat{T}$ and inverse Fourier-Dunkl transform $\widecheck{T}$ are defined as  $ \widehat{T}(f)= T(\widehat{f})$ and $\widecheck{T}(f)=T(\widecheck{f})$ for all $f\in \Sc(\R^n \times \R^d)$, respectively. 
\end{defn}
\begin{defn}\cite{JPSH}
	Let $T\in \Sc'(\R^n\times \R^d)$ be a tempered distribution and $h\in \Sc(\R^n \times \R^d)$. Define the Fourier-Dunkl convolution  $h \ast_{FD} T$ by
	\be (h\ast_{FD} T) (f)=T(\widetilde{h} \ast_{FD}f)\ee for all $f\in \Sc(\R^n \times \R^d)$, where $\widetilde{h}(x, y)=h(-x, -y)$.
\end{defn}
	Note that we can generalize  \eqref{lem:con} and   \eqref{prop:young} to the tempered distributions in the sense that
	$$(f\ast_{FD} T \widehat{)}=\widehat{f}\widehat{T}$$ for all $f\in \Sc(\R^n \times \R^d)$ and
	$$\|f \ast_{FD} T \|_{L^\infty_\K(\R^n\times\R^d)}\leq\|f\|_{L^1_\K(\R^n\times\R^d)}\|T\|_{L^\infty_\K(\R^n\times\R^d)}$$
	whenever $T$ is a tempered distribution given by an $L^\infty$-function which is radial in the second variable.
 
%--------------------------------------------------------------------------------------------------------------------------------------------------------

\section{General theory of restriction theorems}\label{gtrt}
In  this section we recall   crucial ideas to  prove the inequality \eqref{5} and restriction estimate for an orthonormal family of functions (see \eqref{8}).   We note that in     \cite{JPSH} the authors  studied   restriction theory for general surfaces   for  the Fourier-Dunkl transform. However, to make this paper self contained, we give here main ideas and  summery of main results from \cite{JPSH}.

It is easy to verify that $\ecs$ is bounded if and only if $T_S:=\ecs\ecs^*$ is bounded from $L^p$ to $L^{p'},$ where $ \ecs^*$ is the adjoint operator of $ \ecs.$ In \cite{RS} and \cite{ES}, the authors proved the boundedness of $T_S$ (in the case of $\K \equiv 0$) using complex interpolation. In our case, i.e., for any $\K \geq 0$, we also  consider $T_S:= \ecs\ecs^*,$ and prove the boundedness of $T_S$ from $L^p_\K(\R^n\times\R^d)$ into $L^{p'}_\K(\R^n\times\R^d).$ 

In order to prove the boundedness of $T_S,$ we will define an analytic family   of operators $(T_z)$ on the strip $-\la_0\leq \Re(z)\leq0$ with $\la_0>1$ and $T_{-1} = T_S$ and show  that
\begin{equation}\label{7}
	\|T_{is}\|_{L^2_\K(\R^n\times\R^d)\rightarrow L^2_\K(\R^n\times\R^d)}\leq M_0 e^{a|s|}, \|T_{-\la_0+is}\|_{L^1_\K(\R^n\times\R^d)\rightarrow L^\infty_\K(\R^n\times\R^d)}\leq M_1 e^{b|s|},
\end{equation}
for some $\la_0>1$ and $a, b\geq 0$.  Then Stein's complex interpolation theorem implies $T_S = T_{-1}$ is bounded from $L^p_\K(\R^n\times\R^d)$ to $L^{p'}_\K(\R^n\times\R^d)$ with $p = \frac{2\la_0}{\la_0+1}$ and hence we get the inequality \eqref{5} with $p = \frac{2\la_0}{\la_0+1}.$

Note that, as an application of H\"older's inequality, $T_S$ is merely bounded from $L^p_\K(\R^n\times\R^d)$ to $L^{p'}_\K(\R^n\times\R^d),$ if and only if for any $W_1, W_2\in L^{\frac{2p}{2-p}}_\K(\R^n\times\R^d)$, the operator $W_1T_SW_2$ is bounded from $L^2_\K(\R^n\times\R^d)$ to $L^2_\K(\R^n\times\R^d)$ with the estimate
$$\|W_1T_SW_2\|_{L^2_\K\rightarrow L^2_\K} \leq C \|W_1\|_{L^{\frac{2p}{2-p}}_\K} \|W_2\|_{L^{\frac{2p}{2-p}}_\K}.$$
Moreover,  if we have estimates \eqref{7}, then  $W_1T_SW_2$ is more than a mere bounded operator on $L^2_\K,$ namely, it belongs to a Schatten class due to an interpolation idea of Frank-Sabin \cite{FS} in Schatten spaces,  see the following proposition whose proof  can be found in \cite{JPSH}.

% \textcolor{red}{Note that, as an application of H\"older's inequality, to prove $T_S$ is merely bounded from $L^p_\K(\R^n\times\R^d)$ to $L^{p'}_\K(\R^n\times\R^d),$ it is enough to show that for any $W_1, W_2\in L^{\frac{2p}{2-p}}_\K(\R^n\times\R^d)$, the operator $W_1T_SW_2$ is bounded from $L^2_\K(\R^n\times\R^d)$ to $L^2_\K(\R^n\times\R^d)$ with the estimate$$\|W_1T_SW_2\|_{L^2_\K\rightarrow L^2_\K} \leq C \|W_1\|_{L^{\frac{2p}{2-p}}_\K} \|W_2\|_{L^{\frac{2p}{2-p}}_\K}.$$ Moreover,  if we have estimates \eqref{7},  a duality principle shows that $W_1T_SW_2$ is more than a mere bounded operator on $L^2_\K,$ namely, it belongs to a Schatten class.}

%\textcolor{red}{Using the idea of  Frank-Sabin \cite{FS}, in order to prove the Restriction Problem  (\ref{5}) and  restriction estimate for an orthonormal family of functions in the context of   the Fourier-Dunkl transform,  the following proposition and lemma play an important role.}
\begin{prop}\label{A}
	Let $(T_z)$ be an analytic family of operators on $\R^n \times \R^d$ in the sense of Stein defined on the strip $-\la_0 \leq \text{Re}(z) \leq 0$ for some $\la_0 > 1.$ Assume that we have the following bounds
	\begin{equation}\label{B1}
		\|T_{is}\|_{L^2_\K(\R^n\times\R^d)\rightarrow L^2_\K(\R^n\times\R^d)}\leq M_0 e^{a|s|},\quad  \|T_{-\la_0+is}\|_{L^1_\K(\R^n\times\R^d)\rightarrow L^\infty_\K(\R^n\times\R^d)}\leq M_1 e^{b|s|},
	\end{equation}
	for all $s\in \R,$ for some $a, b, M_0, M_1 \geq 0.$ Then, for all $W_1, W_2 \in L_\K^{2\la_0}(\R^n\times \R^d)$,  the operator $W_1T_{-1}W_2$ belongs to $\s^{2\la_0}(L^2_\K(\R^n\times\R^d))$ and we have the estimate
	\begin{equation}\label{B2}
		\|W_1T_{-1}W_2\|_{\s^{2\la_0}(L^2_\K(\R^n\times\R^d))}\leq M^{1-\frac{1}{\la_0}}_0M^\frac{1}{\la_0}_1 \|W_1\|_{L^{2\la_0}_{\K}(\R^n\times\R^d)} \|W_2\|_{L^{2\la_0}_{\K}(\R^n\times\R^d)}.
	\end{equation}
\end{prop}
If we have the Schatten norm estimates for the operator $W_1T_SW_2$, then the following lemma gives the restriction estimates for the system of orthonormal functions.
\begin{lem}[Duality Principle, \cite{JPSH}]\label{B} 
	Let $\Hc$ be a Hilbert space. Assume that $A$ be a bounded operator from $\Hc$ to $L^{p'}_{\K}(\R^n \times \R^d)$ for some $1\leq p\leq 2$ and let $\ap\geq 1$. Then the following statements are equivalent.
	\begin{enumerate}
		\item There is a constant $C > 0$ such that
		\begin{equation}\label{DI}
			\|WAA^*\overline{W}\|_{\s^\ap\left(L^2_{\K}(\R^n\times \R^d)\right)} \leq C \|W\|^2_{L^\frac{2p}{2-p}_\K(\R^n \times \R^d)},
		\end{equation}
		for all $W\in {L^\frac{2p}{2-p}_\K(\R^n \times \R^d)}.$ %where the function $W$ is interpreted as an operator which acts by multiplication.\\
		\item  There is a constant $C'>0$ such that for any orthonormal system $(f_j)_{j\in J}$ in $\Hc$ and any sequence $(n_j)_{j\in J}\subset \C,$ we have
		\begin{equation}\label{DI2}
			\left\|\sum_{j\in J}n_j\left|Af_j\right|^2\right\|_{L^{\frac{p'}{2}}_{\K}(\R^n \times \R^d)}\leq C' \left(\sum_{j\in J}|n_j|^{\ap'}\right)^{1/\ap'}.
		\end{equation}
		Moreover, the values of the optimal constants $C$ and $C'$ coincide.
	\end{enumerate} 
\end{lem}
 
As we discussed earlier, once we have estimates \eqref{7}, Proposition \ref{A} implies  that $W_1T_SW_2 $ belongs to a Schatten class,  and further using Lemma \ref{B}, we will get the following restriction estimate for an orthonormal family of functions
\begin{equation}\label{8}
	\left\|\sum_jn_j|\ecs f_j|^2\right\|_{L^r_\K(\R^n\times\R^d)}\leq C \left(\sum_j|n_j|^\bt\right)^{1/\bt},
\end{equation}
where $\{f_j\}$ is a sequence of orthonormal function in $L^2(S, h^2_\K(y)d\sigma), r = \frac{\la_0}{\la_0-1}~~\text{and}~~\bt=\frac{2\la_0}{2\la_0-1}.$

Thus  the  main goal is to provide a sufficient condition for the estimates \eqref{7} to hold for the the general surface.  Consider the quadratic surfaces  
\begin{equation}\label{9}
	\mbox{Let}~~ S= S_r=\{(x, y) \in \R^n\times \R^d: P(x, y)=r\},
\end{equation}
where $P(x, y)$ is a polynomial of degree two with real coefficients and $r$ is a real constant. We assume that $P$ is not a function of fewer than $n$ variables so that $S$ is a $(n+d-1)$-dimensional $C^\infty$-manifold with the canonical measure $d\mu$ associated to the function $P$ given by
\begin{equation}\label{10}
	d\mu_r(x, y)=\frac{dx_1 dx_2 \cdots dx_{n-1} dy}{|\partial P/\partial x_n|}
\end{equation}
in any neighborhood in which $\frac{\partial P}{\partial x_n}\neq 0$ so that $S$ may be described by giving $x_n$ as a function of $x_1, x_2, \ldots, x_{n-1}, y_1, y_2, \ldots, y_{d}.$

%\textcolor{red}{Now, as an application of H\"older's inequality, to prove $T_S$ is merely bounded from $L^p_\K(\R^n\times\R^d)$ to $L^{p'}_\K(\R^n\times\R^d),$ it is enough to show that for any $W_1, W_2\in L^{\frac{2p}{2-p}}_\K(\R^n\times\R^d)$, the operator $W_1T_SW_2$ is bounded from $L^2_\K(\R^n\times\R^d)$ to $L^2_\K(\R^n\times\R^d)$ with the estimate	$$\|W_1T_SW_2\|_{L^2_\K\rightarrow L^2_\K} \leq C \|W_1\|_{L^{\frac{2p}{2-p}}_\K} \|W_2\|_{L^{\frac{2p}{2-p}}_\K}.$$
	%	But if we have estimates \eqref{7}, Proposition \ref{A} shows that $W_1T_SW_2$ is more than a mere bounded operator on $L^2_\K,$ namely that it belongs to a Schatten class, as a result (using Lemma \ref{B}) we will get the following restriction estimate for an orthonormal family of functions
	%	\begin{equation}\label{8}
		%		\left\|\sum_jn_j|\ecs f_j|^2\right\|_{L^r_\K(\R^n\times\R^d)}\leq C \left(\sum_j|n_j|^\bt\right)^{1/\bt},
		%	\end{equation}
	%	where $\{f_j\}$ is a sequence of orthonormal function in $L^2(S, h^2_\K(y)d\sigma), r = \frac{\la_0}{\la_0-1}~~\text{and}~~\bt=\frac{2\la_0}{2\la_0-1}.$
	%	As we have said earlier that we prove a sufficient condition for the estimates \eqref{7} to hold for the cone surface, and hence we get the following:}

As we have discussed  earlier, to prove inequality \eqref{5}, it is enough to prove $T_S=\ecs\ecs^*$ is bounded from  $L^p_\K(\R^n\times\R^d)$ to $L^{p'}_\K(\R^n\times\R^d).$ In order to do that, we define an analytic family of operators. We consider the generalized function $$G_z(x, y) = w(z)(P(x, y)-r)
^z_+,$$ where $w(z)$ is an appropriate analytic function with a simple zero at $z=-1$ and
\begin{align}\label{generalized}
	s^z_+=\begin{cases}s^z, ~~&\mbox{if}~s>0,\\ 0,~~&\mbox{if}~s\leq 0. \end{cases}
\end{align}%  From the definition of $d\mu_r$ given in \eqref{10} and $G_z,$ we have $$\int_{\R^n\times\R^d}G_z(x, y)\phi(x, y)dxh^2_\K(y)dy=w(z)\int_\R\left(\int_{S_t}\phi(x, y) h^2_\K(y) d\mu_t(x, y)\right)(t-r)^z_+ dt,$$for any Schwartz class functions $\phi$ on $\R^n\times\R^d.$ 
Further, using the one-dimensional analysis of $(t-r)^z_+,$ (see \cite{GS}), we have
$$\lim_{z\rightarrow -1}\int_{\R^n\times\R^d}G_z(x, y)\phi(x, y)dxh^2_\K(y)dy = C \int_{S_r}\phi(x, y) h^2_\K(y) d\mu_r(x, y),$$
for all $\phi\in \Sc(\R^n\times\R^d)$.  Keeping this observation in mind, we consider the analytic family of operators $T_z g:= (G_z \widehat{g})^\vee,$ for all $g\in \Sc(\R^n\times\R^d)$, i.e.,
$$T_z g(x, y)= \frac{1}{c_\K(2\pi)^{n/2}} \int_{\R^n\times\R^d}G_z(\xi, \zeta) \widehat{g}(\xi, \zeta)e^{ix\cdot\xi}E_\K(iy, \zeta) h^2_\K(\zeta)d\xi d\zeta.$$
Then  $T_{-1} =C T_S.$ Since $G_{is}$ is bounded, we have that %$$\|T_{is}g\|_{L^2_\K(\R^n\times\R^d)}=\|G_{is}\widehat{g}\|_{L^2_\K(\R^n\times\R^d)}\leq|w(is)| \|g\|_{L^2_\K(\R^n\times\R^d)},$$ i.e., 
\begin{equation}\label{11}
	\|T_{is}\|_{L^2_\K(\R^n\times\R^d)\rightarrow L^2_\K(\R^n\times\R^d)} \leq |w(is)|,
\end{equation}
where we  choose $w(z)$ so that $|w(is)|$ has at most exponential growth at infinity in $s$. Again, if we have 
\begin{equation}\label{12}
	\|T_{-\la_0+is}\|_{L^1_\K(\R^n\times\R^d)\rightarrow L^\infty_\K(\R^n\times\R^d)}\leq C e^{b|s|}
\end{equation}
for some $\la_0> 1,$ then Stein's complex interpolation gives that the boundedness of $T_S = T_{-1}$ from $L^p_\K(\R^n\times\R^d)\rightarrow L^{p'}_\K(\R^n\times\R^d)$ for $p = \frac{2\la_0}{\la_0+1}.$ Then using Proposition \ref{A} and Lemma \ref{B}, we   get inequality \eqref{8}.
%\subsection{A sufficient condition for \eqref{12} to hold}
%If we have to follow very closely the article \cite{RS} then we have to
%\begin{enumerate}
%\item define  some kind of convolution call it Fourier-Dunkl Convolution $(f\ast_{FD}g)$ of functions on $\R^n\times\R^d$ satisfying $\widehat{(f\ast_{FD}g)}= \hat{f}\hat{g},$ where $\hat{f}$ and $\hat{g}$ are the Fourier-Dunkl transform of $f~~ \text{and}~~ g.$ Here we are proposing a definition of Fourier-Dunkl convolution $$(f\ast_{FD}g)(x, y) = \int_{\R^n\times\R^d}\tau_yf(x-\xi, -\zeta)g(\xi, \zeta)h^2_\K(\zeta)d\xi d\zeta,$$  where $\tau_y f(\xi, \zeta)$ is the Dunkl transform of $f$ in second variable $y.$
%We have to verify that $\widehat{(f\ast_{FD} g)} = \hat{f}\hat{g}.$
%\item extend Fourier-Dunkl transform to the tempered distributions $\phi$ on $\Sc(\R^n\times\R^d)$ and Fourier-Dunkl convolution of a tempered distribution $\phi$ and Schwartz class function $f$ such that $$(\phi g\check) = \widehat{\phi}\ast_{FD} \widehat{g},$$
%whenever the distribution $\phi$ is given by a function $\phi(x, y).$
%\item have $\|\phi\ast_{FD} g\|_{L^\infty_\K(\R^n\times\R^d)} \leq \|\phi\|_{L^\infty_\K(\R^n\times\R^d)}\|g\|_{L^1_\K(\R^n\times\R^d)}$.
%\end{enumerate}
%Suppose we have the above $(1)$ to $(3).$
From the definition of $T_z$ we see that
$$T_z g= (G_z \widehat g\widecheck{)}=   g \ast_{FD} \widecheck{G}_z $$
and hence
\begin{equation}\label{13}
	\|T_z g\|_{L^\infty_\K(\R^n\times\R^d)}\leq\| \widecheck{G}_z \|_{L^\infty_\K(\R^n\times\R^d)} \|g\|_{L^1_\K(\R^n\times\R^d)} ,
\end{equation}
provided $G_z$ is radial in second variable. In view of \eqref{13}, we have the following result for general surfaces.
\begin{theorem}\label{C}
	Suppose  that the generalized function $G_z$ is radial in second variable, $\widecheck{G}_z$  is bounded for $\Re(z) = -\la_0$ with $\la_0>1$ and $\| \widecheck{G}_{-\la_0+is}\|_{L^\infty}$ has at most exponential growth at infinity with respect to $s$. Then \eqref{5} and \eqref{8} holds for $p=\frac{2\la_0}{\la_0+1}, r=\frac{\la_0}{\la_0-1}$ and $\bt= \frac{2\la_0}{2\la_0-1}.$
\end{theorem}

%\begin{rem}Using these results for special surfaces $S$, we can also prove (we believe) the  Strichartz inequalities for the following cases:
%\begin{enumerate}
%\item Dunkl Laplacian (Schr\"ondinger's equation).
%\item Square root of the Dunkl Laplacian (Wave equation)
%\item Pseudo-relativistic (The Klein–Gordon equation)
%
%\end{enumerate}
%\end{rem}

%-------------------------------------------------------------------------------------------------------------------------------------------------

\section{Fourier-Dunkl transform of certain generalized functions}\label{sec3}
In this section, we will find the Fourier-Dunkl transform of certain generalized functions, namely $(x_n-|x'|^2+|y|^2)_+^z,~~(1-|x|^2-|y|^2)_+^z$, and $(1-|x|^2+|y|^2)_+^z,$ where $x'= (x_1, \ldots, x_{n-1}).$

\subsection{Fourier-Dunkl transform of $(x_n-|x'|^2+|y|^2)_+^z$}
Let $P'(x', y)= |x'|^2-|y|^2.$ Using the fact that (see \cite{GS}, p. 360)
$$\int_0^\infty t^z e^{-itw} dt= \Gm(z+1)ie^{iz\frac{\pi}{2}}(-w+i0)^{-z-1},$$
 we obtain
\begin{equation}\label{fd1}
\int_{-\infty}^\infty e^{-ix_n\cdot \xi_n}(x_n-P'(x', y))_+^z dx_n = \Gm(z+1)ie^{iz\frac{\pi}{2}} e^{-i\xi_nP'(x', y)}(-\xi_n+i0)^{-z-1}.
\end{equation}
Moreover, it can be shown that 
\begin{equation}\label{fd2}
\frac{1}{\sqrt{2\pi}}\int_{-\infty}^\infty e^{itx^2} e^{ixy} dx =\sqrt{2} |t|^{-\frac{1}{2}} e^{\frac{i\pi}{4} \sgn(t)} e^{-\frac{iy^2}{4t}},
\end{equation}
the Fourier transform is taken in the sense of distribution. Also from   \cite{MR}, it is known that 
$$\int_{\R^d} e^{-a|y|^2} E_\K(-i\z, y) h^2_\K(y) dy= \frac{c_\K}{2^{\frac{d+2\gk}{2}}}|a|^{-\frac{d+2\gk}{2}}e^{-\frac{d+2\gk}{2}i Arg(a)}e^{-\frac{1}{4a}|\z|^2},$$
whenever $\Re(a)>0.$
Let $a=r+it$ and  leeting $r\to 0^+$ in the sense of distribution,   we get 
\begin{equation}\label{fd3}
\int_{\R^d} e^{-it|y|^2}E_\K(-i\z, y) h^2_\K(y) dy= \frac{c_\K}{2^{\frac{d+2\gk}{2}}} |t|^{-\frac{d+2\gk}{2}} e^{-\frac{d+2\gk}{4}\pi i \sgn(t)} e^{\frac{i}{4t}|\z|^2}.
\end{equation}  
We are now ready to compute the Fourier-Dunkl transform of $(x_n-P')_+^z.$ Using the definition (\ref{eq:Dunkltr}), we get 
%\begin{align*} 
$$(x_n-P'\widehat{)^z_+}(\xi, \z)= \frac{1}{c_\K(2\pi)^\frac{n}{2}}\int_{\R^d}\int_{\R^{n-1}} \Big(\int_\R(x_n-P')_+^z e^{-ix_n\cdot \xi_n}dx_n\Big) e^{-ix'\cdot\xi'}E_\K(-i\z, y) dx' h^2_\K(y) dy,$$
 where $\xi'= (\xi_1, \ldots, \xi_{n-1}).$ From \eqref{fd1}, we get
%\end{align*}
\begin{multline*}
\frac{(x_n-P'\widehat{)^z_+}(\xi, \z)}{\Gm(z+1)}= \frac{ie^{iz\frac{\pi}{2}}(-\xi_n+i0)^{-z-1}}{c_\K(2\pi)^\frac{n}{2}}\int_{\R^d}\left(\int_{\R^{n-1}}  e^{-i\xi_nP'(x', y)} e^{-ix'\cdot\xi'} dx'\right) E_\K(-i\z, y) h^2_\K(y) dy.
\end{multline*}
Using \eqref{fd2}, the above equals to 
\begin{align*}
\frac{(x_n-P'\widehat{)^z_+}(\xi, \z)}{\Gm(z+1)}&= \frac{ie^{iz\frac{\pi}{2}}}{c_\K\sqrt{2\pi}}(-\xi_n+i0)^{-z-1} e^{-\frac{i(n-1)\pi}{4}\sgn{\xi_n}} e^{\frac{i|\xi'|^2}{4\xi_n}}\left(\frac{|\xi_n|}{2}\right)^{-\frac{n-1}{2}}\\ &\qquad \times\int_{\R^d}e^{i\xi_n|y|^2}E_\K(-i\z, y)h^2_\K(y) dy.
\end{align*}
Combining this with \eqref{fd3}, we obtain
\begin{align}\label{DFgp1}
\frac{(x_n-P'\widehat{)^z_+}(\xi, \z)}{\Gm(z+1)}&= \frac{ie^{iz\frac{\pi}{2}}}{c_\K\sqrt{\pi}2^\frac{d+2\gk-n+2}{2}} e^{\frac{i\pi}{4}(d+2\gk-n+1)\sgn \xi_n } e^{\frac{i}{4\xi_n}P'(\xi', \z)}|\xi_n|^{-\frac{(n+d+2\gk-1)}{2}}(-\xi_n+i0)^{-z-1},
\end{align}
where $P'(\xi', \z)= |\xi'|^2-|\z|^2.$ This is  bounded if and only if $\Re(z)= - \frac{(d+n+2\gk+1)}{2}$ and the growth as Im$(z)\to 0$ is exponential.   

For computing Fourier-Dunkl transform of $(1-|x|^2-|y|^2)_+^z$ and $(1-|x|^2+|y|^2)_+^z,$ we use the method of \cite{GS} which is explained in the following subsequent subsections. 

\subsection {Fourier-Dunkl transform of $(c^2+P+i0)^z$ and $(c^2+P-i0)^z$} 
In this subsection, we first define  the generalized functions $(c^2+P+i0)^z$ and $(c^2+P-i0)^z$ and then we compute their Fourier-Dunkl transforms using the method of \cite{GS} for the class of $P$ which are of the form 
\begin{equation}\label{qf1}P(x, y)=\sum_{j=1}^n \alpha_j x_j^2+\bt |y|^2,\end{equation} 
where $(x, y) \in \R^n \times \R^d, \alpha_j\in \R, j=1, \ldots ,n,$ and $\bt \in \R$. For any real quadratic form $P$, we define the generalized functions $(c^2+P+i0)^z$ and $(c^2+P-i0)^z$ by
\begin{equation}\label{limge1}
(c^2+P+i0)^z= \lim_{\ve\to 0}(c^2+P+i\ve P_0)^z
\end{equation}
and
\begin{equation}\label{limge2}
(c^2+P-i0)^z= \lim_{\ve\to 0}(c^2+P-i\ve P_0)^z,
\end{equation}
where $\ve>0$ and $P_0$ is a positive definite quadratic form.  The existence of the limits in \eqref{limge1} and \eqref{limge2} for $c=0$  (see \cite[Section 2.4]{GS}) and the existence for $c\neq 0$ follows from the absence of singular points on the $c^2+P=0$ hyper surface. From  \cite[Subsection 2.4]{GS}, we have that 
\begin{equation}\label{ar1}
(c^2+P)^\la_+= \frac{i}{2 \sin \la\pi} \{e^{-i\la\pi}(c^2+P+i0)^\la-e^{i\la\pi}(c^2+P-i0)^\la\},
\end{equation}
i.e., the Fourier-Dunkl transform of $(c^2+P)^\la_+$  is   the linear combination of that of $(c^2+P+i0)^\la$ and $(c^2+P-i0)^\la$. The required generalized functions, namely $(1-|x|^2\pm |y|^2)^z_+$ are of the form $(c^2+P)^z_+$.   Thus we will first find the  Fourier-Dunkl transforms of $(c^2+P+i0)^\la$ and $(c^2+P-i0)^\la$.

We start by considering the generalized function $(c^2+P)^z$ for a quadratic form $P$ as in \eqref{qf1} with $\ap_j>0, j=1, \ldots, n$ and $\bt>0$. The Fourier-Dunkl transform of this generalized function for $\Re(z)<-\frac{N_\K}{2}$ is given by
\begin{equation}
(c^2+P\widehat{)^z}(\xi, \z)=\frac{1}{c_\K(2\pi)^{n/2}}\int_{\R^n}\int_{\R^d} (c^2+P)^z E_\K(-i\zeta, y)  e^{-ix\cdot \xi} h^2_\K(y) dydx,~~~(\xi, \z)\in \R^n\times \R^d.
\end{equation}
A simple change of variables gives that
\begin{align}
(c^2+P\widehat{)^z}(\xi, \z)&=\frac{1}{c_\K(2\pi)^{n/2}\sqrt{D}}\int_{\R^n}\int_{\R^d} (c^2+|x|^2+|y|^2)^z E_\K(-i\z/\sqrt{\bt}, y)  e^{-ix\cdot \xi_a} h^2_\K(y) dydx\nonumber \\
&=\frac{1}{\sqrt{D}}(c^2+|x|^2+|y|^2\widehat{)^z}(\xi_a, \z/\sqrt{\bt}),\label{fdgf1}
\end{align}
where $\xi_a=\left(\frac{\xi_1}{\sqrt{a_1}}, \cdots , \frac{\xi_n}{\sqrt{a_n}}\right)$, $D=a_1\cdots a_n \bt^{d+2\gk}$, and $\Re(z)<-\frac{N_\K}{2}$. Therefore,  by writing the integrals over $\R^n$ and $\R^d$ in polar coordinates, we get
\begin{align*}
(c^2+|x|^2+|y|^2\widehat{)^z}(\xi, \z)&= \frac{1}{c_\K(2\pi)^n/2}\int_0^\infty\int_0^\infty(c^2+r^2_1+r^2_2)^z\left(\int_{S^{n-1}}e^{-ir_1x'\cdot\xi}d\si(x')\right)\\ &\qquad\qquad\times\left(\int_{S^{d-1}}E_\K(-i\z, r_2y')h^2_\K(y')d\si(y')\right)r_1^{n-1}r_2^{d+2\gk-1}dr_1dr_2.
\end{align*}  Here  $S^{n-1}$ and  $S^{d-1}$  denotes  the unit spheres respectively in $\R^{n-1}$ and $\R^{d-1}$. Using the formula (see \cite[Proposition 2.8]{CFD5})
\begin{equation}
\int_{S^{d-1}}E_\K(-\z, y)h_\K^2(y) d\sigma(y)=\sigma_{d-1, \K}\Gm\left(\frac{d+2\gk}{2}\right)\left(\frac{|\z|}{2}\right)^{-\frac{d+2\gk}{2}+1}J_{\frac{d+2\gk}{2}-1}(|\z|)
\end{equation}
where $J_\nu(r)$ is the Bessel function of the first kind and $\displaystyle \si_{d-1, \K}= \int_{S^{d-1}}h^2_\K(y) d\si(y)$, we get that 
\begin{align*}
(c^2+|x|^2+|y|^2\widehat{)^z}(\xi, \z)&=C_{n,d,\K} \Big(\frac{|\xi|}{2}\Big)^{-\frac{n}{2}+1}\Big(\frac{|\z|}{2}\Big)^{-\frac{d+2\gk}{2}+1} \int_0^\infty\int_0^\infty (c^2+r^2_1+r^2_2)^z r_1^{\frac n 2}r_2^{\frac{d+2\gk}{2}}\\& \qquad \times J_{\frac{n}{2}-1}(r_1|\xi|)J_{\frac{d+2\gk}{2}-1}(r_2|\z|)dr_1dr_2,
\end{align*}
where $\displaystyle C_{n, d, \K}= \frac{\si_{n-1}\si_{d-1, \K}\Gm\big(\frac{n}{2}\big)\Gm\big(\frac{2\gk+d}{2}\big)}{c_\K (2\pi)^{n/2}}$ with $ \si_{n-1}=\si_{n-1, 0}.$ Using the change variable $r_1=r\cos \theta, r_2=r \sin \theta$ and then using Sonine's second finite integral formula (see \cite[p. 376 ]{WT}) for the above integral, we obtain that
$$
(c^2+|x|^2+|y|^2\widehat{)^z}(\xi, \z)= \frac{C_{n, d, \K}}{|s|^{(\frac{N_\K}{2}-1)}}\int_0^\infty (c^2+r^2)^z r^\frac{N_\K}{2}J_{\frac{N_\K}{2}-1}(r|s|)dr,$$
where $s^2=|\xi|^2+|\z|^2$.  Using a formula from \cite[p. 434 ]{WT}, we further get   
\begin{equation}\label{wkfb}
(c^2+|x|^2+|y|^2\widehat{)^z}(\xi, \z)=\frac{C_{n, d, \K}2^{z+1}}{\Gm(-z)}\left(\frac{c}{|s|}\right)^{\frac{N_\K}{2}+z} K_{\frac{N_\K}{2}+z}(c|s|),
\end{equation}
for $\Re(z)<-N_\K/2$, where $$K_\la(z)=\frac{\pi}{2\sin \la \pi}\left[e^{\frac 1 2 i\la \pi}J_{-\la}(iz)-e^{-\frac 1 2 i\la \pi}J_{\la}(iz)\right].$$ For other values of $z$ in equation \eqref{wkfb} remains valid by analytic continuation in $z.$ In view of \eqref{fdgf1}, we get the Fourier-Dunkl transform of $(c^2+P)^z$ as follows:
\begin{equation}\label{equ:conz}
(c^2+P\widehat{)^z}(\xi, \z)= \frac{C_{n, d, \K}2^{z+1}}{\sqrt{D}\Gm(-z)}\Big(\frac{c}{Q^\frac{1}{2}(\xi, \z)}\Big)^{\frac{N_\K}{2}+z}K_{\frac{N_\K}{2}+z}(cQ^\frac{1}{2}(\xi, \z)),
\end{equation}
where $\displaystyle Q(\xi, \z)= \sum_{i=1}^n\frac{\xi_i^2}{\ap_i}+\frac{1}{\bt}|\z|^2.$

Now let $P$ be any real quadratic form as in \eqref{qf1}.  Let $(c^2+P+i0)^z$ and $(c^2+P-i0)^z$ be defined by \eqref{limge1} and \eqref{limge2}, respectively. If the quadratic form $\mathscr{P}$ lies in the ``upper half-plane,'' then its dual $\mathscr{Q}$ lies in the ``lower half-plane'' with the same arguments given in \cite[Subsection 2.8]{GS}. Therefore, according to the uniqueness of analytic continuation, equation \eqref{equ:conz} implies that 
\begin{equation}\label{fdgf2}
(c^2+P+i0\widehat{)^z}(\xi, \z)= \frac{C_{n, d, \K}2^{z+1}c^{\frac{N_\K}{2}+z}}{\sqrt{D}\Gm(-z)}\;\frac{K_{\frac{N_\K}{2}+z}\big(c(Q(\xi, \z)-i0)^\frac{1}{2}\big)}{(Q(\xi, \z)-i0)^{\frac{1}{2}{\big(\frac{N_\K}{2}+z}\big)}}.\end{equation}
Here by $\sqrt{D}$ we denote the analytic continuation from the sheet on which this function is positive for $\ap_j>0, j=1, \ldots, n$ and $\bt>0$. Similarly, it can be shown that 
\begin{equation}\label{fdgf3}
(c^2+P-i0\widehat{)^z}(\xi, \z)= \frac{C_{n, d, \K}2^{z+1}c^{\frac{N_\K}{2}+z}}{\sqrt{D}\Gm(-z)}\;\frac{K_{\frac{N_\K}{2}+z}\big(c(Q(\xi, \z)+i0)^\frac{1}{2}\big)}{(Q(\xi, \z)+i0)^{\frac{1}{2}{\big(\frac{N_\K}{2}+z}\big)}}.
\end{equation}

%-------------------------------------------------------------------------------------------------------------------------------------------------------------
\subsection{Fourier-Dunkl transform of $(1-|x|^2-|y|^2)_+^z$}

Let $\ap_j=-1,~~\bt=-1$ for all $j=1, \ldots, n$ and let $P_1=-|x|^2-|y|^2$ and $Q_1(\xi, \z)= -|\xi|^2-|\z|^2.$ Then $\sqrt{D}= \sqrt{|D|}e^{\frac{1}{2}N_\K i\pi}$ and the formula \eqref{fdgf2} gives
\begin{equation}\label{fd11}
(c^2+P_1+i0\widehat{)^z}(\xi, \z)= \frac{C_{n, d, \K} 2^{z+1}c^{\frac{N_\K}{2}+z}}{\sqrt{|D|}e^{\frac{1}{2}N_\K i\pi}\Gm(-z)}\;\frac{K_{\frac{N_\K}{2}+z}(c(Q_1(\xi, \z)-i0)^\frac{1}{2})}{(Q_1(\xi, \z)-i0)^{\frac{1}{2}\big(\frac{N_\K}{2}+z\big)}}.
\end{equation}
Similarly, the formula \eqref{fdgf3} gives
\begin{equation}\label{fd12}
(c^2+P_1-i0\widehat{)^z}(\xi, \z)=  \frac{C_{n, d, \K} 2^{z+1}c^{\frac{N_\K}{2}+z}}{\sqrt{|D|}e^{-\frac{1}{2}N_\K i\pi}\Gm(-z)}\frac{K_{\frac{N_\K}{2}+z}(c(Q_1(\xi, \z)+i0)^\frac{1}{2})}{(Q_1(\xi, \z)+i0)^{\frac{1}{2}\big(\frac{N_\K}{2}+z\big)}}.  
\end{equation}
%where $\sqrt{|\D|}e^{-\frac{1}{2}N_\K i\pi}=e^{-\frac{1}{2}N_\K i\pi}.$ 
Using the identity
\begin{equation}\label{gi}
\begin{cases}
(x+i0)^\la&= x^\la_+ +e^{i\la\pi}x^\la_-,\\
(x-i0)^\la&= x^\la_+ +e^{-i\la\pi}x^\la_-
\end{cases}
\end{equation}
and expanding $K_\la$ in power series we get that
\begin{equation}\label{kpe}
\frac{K_\la(c(x+i0)^\frac{1}{2})}{(x+i0)^\frac{\la}{2}}=\frac{ K_\la(cx_+^\frac{1}{2})}{x_+^\frac{\la}{2}}+\frac{\pi}{2\sin \la\pi}\left[e^{-i\la\pi}\frac{J_{-\la}(cx^\frac{1}{2}_-)}{x_-^\frac{\la}{2}}-\frac{J_{\la}(cx^\frac{1}{2}_-)}{x_-^\frac{\la}{2}}\right]
\end{equation}
and
\begin{equation}\label{kpe2}
  \frac{K_\la(c(x-i0)^\frac{1}{2})}{(x-i0)^\frac{\la}{2}}=\frac{ K_\la(cx_+^\frac{1}{2})}{x_+^\frac{\la}{2}}+\frac{\pi}{2\sin \la\pi}\left[e^{i\la\pi}\frac{J_{-\la}(cx^\frac{1}{2}_-)}{x_-^\frac{\la}{2}}-\frac{J_{\la}(cx^\frac{1}{2}_-)}{x_-^\frac{\la}{2}}\right]. 
\end{equation}
Let us now express \eqref{fd11} and \eqref{fd12} in terms of ${Q_1}_+(\xi, \z)$ and ${Q_1}_-(\xi, \z)$ by using the expressions of \eqref{kpe} and \eqref{kpe2} as follows:
\begin{multline}\label{ffdt}
(c^2+P_1+i0\widehat{)^z}(\xi, \z)=\frac{C_{n, d, \K} 2^{z+1}c^{\frac{N_\K}{2}+z}}{\sqrt{|D|}e^{\frac{1}{2}N_\K i\pi}\Gm(-z)} \left[\frac{K_{\frac{N_\K}{2}+z}(c{Q_1}^\frac{1}{2}_+(\xi, \z))}{{Q_1}^{\frac{1}{2}(\frac{N_\K}{2}+z)}_+(\xi, \z)}\right.\\+\left.\frac{\pi}{2\sin (\frac{N_\K}{2}+z)\pi}\left\{e^{i(\frac{N_\K}{2}+z)\pi}\frac{J_{-(\frac{N_\K}{2}+z)}(c{Q_1}^\frac{1}{2}_- (\xi, \z))}{{Q_1}^{\frac{1}{2}(\frac{N_\K}{2}+z)}_-(\xi, \z)}-\frac{J_{(\frac{N_\K}{2}+z)}(c{Q_1}^\frac{1}{2}_-(\xi, \z))}{{Q_1}^{\frac{1}{2}(\frac{N_\K}{2}+z)}_-(\xi, \z)}\right\}\right],
\end{multline}
and
\begin{multline}\label{ffdt2}
(c^2+P_1-i0\widehat{)^z}(\xi, \z)=\frac{C_{n, d, \K} 2^{z+1}c^{\frac{N_\K}{2}+z}}{\sqrt{|D|}e^{-\frac{1}{2}N_\K i\pi}\Gm(-z)}   \left[\frac{K_{\frac{N_\K}{2}+z}(c{Q_1}^\frac{1}{2}_+(\xi, \z))}{{Q_1}^{\frac{1}{2}(\frac{N_\K}{2}+z)}_+(\xi, \z)}\right.\\+\left.\frac{\pi}{2\sin (\frac{N_\K}{2}+z)\pi}\left\{e^{-i(\frac{N_\K}{2}+z)\pi}\frac{J_{-(\frac{N_\K}{2}+z)}(c{Q_1}^\frac{1}{2}_-(\xi, \z))}{{Q_1}^{\frac{1}{2}(\frac{N_\K}{2}+z)}_-(\xi, \z)}-\frac{J_{(\frac{N_\K}{2}+z)}(c{Q_1}^\frac{1}{2}_-(\xi, \z))}{{Q_1}^{\frac{1}{2}(\frac{N_\K}{2}+z)}_-(\xi, \z)}\right\}\right].
\end{multline}
Now using the identity \eqref{ar1}, we get
%\begin{align*}
%(c^2+p_1\widehat{)^z}_+(\xi, \z)&=\frac{i}{2 \sin z\pi} \Big[e^{-iz\pi}(c^2+p_1+i0\widehat{)^z}-e^{iz\pi}(c^2+p_1-i0\widehat{)^z}\Big]\\
%&=\frac{i c_{d, \K} 2^{z+1}c^{\frac{N_\K}{2}+z}}{2 \sin z\pi(2\pi)^{\frac{n}{2}}\sqrt{|\D|}\Gm(-z)}\left\{\frac{K_{\frac{N_\K}{2}+z}(c{Q_1}^\frac{1}{2}_+)}{{Q_1}^{\frac{1}{2}(\frac{N_\K}{2}+z)}_+}\big[e^{-iz\pi}e^{-\frac{1}{2}N_\K i\pi}-e^{iz\pi}e^{\frac{1}{2}N_\K i\pi}\big]\right. \\& + \frac{\pi}{2\sin (\frac{N_\K}{2}+z)\pi}\left(\frac{J_{-(\frac{N_\K}{2}+z)}(c{Q_1}^\frac{1}{2}_-)}{{Q_1}^{\frac{1}{2}(\frac{N_\K}{2}+z)}_-}\big[e^{-iz\pi}e^{i(\frac{N_\K}{2}+z)\pi}e^{-\frac{1}{2}N_\K i\pi}-e^{iz\pi}e^{-i(\frac{N_\K}{2}+z)\pi}e^{\frac{1}{2}N_\K i\pi}\big]\right.\\&-\left.\left. \frac{J_{\frac{N_\K}{2}+z}(c{Q_1}^\frac{1}{2}_-)}{{Q_1}^{\frac{1}{2}(\frac{N_\K}{2}+z)}_-}\big[e^{-iz\pi}e^{-\frac{1}{2}N_\K i\pi}-e^{iz\pi}e^{\frac{1}{2}N_\K i\pi}\big]\right)\right\}
%\end{align*}
\begin{equation}\label{lcfdt}
(c^2+P_1\widehat{)^z_+}(\xi, \z)=\frac{i}{2 \sin z\pi} \Big[e^{-iz\pi}(c^2+P_1+i0\widehat{)^z} (\xi, \z)-e^{iz\pi}(c^2+P_1-i0\widehat{)^z}(\xi, \z)\Big].
\end{equation}
By inserting equations \eqref{ffdt} and \eqref{ffdt2} into (\ref{lcfdt}), we obtain 
%\begin{align*}
\begin{align*}(c^2+P_1\widehat{)^z_+}(\xi, \z)&=\frac{i C_{n, d, \K} 2^{z+1}c^{\frac{N_\K}{2}+z}}{2 \sqrt{|D|}\Gm(-z) \sin z\pi}
\left[\frac{K_{\frac{N_\K}{2}+z}(c{Q_1}^\frac{1}{2}_+(\xi, \z))}{{Q_1}^{\frac{1}{2}(\frac{N_\K}{2}+z)}_+(\xi, \z)}\left(-2i \sin\Big(\frac{N_\K}{2}+z\Big)\right)\right. \\&\qquad \qquad \qquad \qquad \qquad \qquad \qquad \qquad \quad + \left. i\pi \frac{J_{\frac{N_\K}{2}+z}(c{Q_1}^\frac{1}{2}_-(\xi, \z))}{{Q_1}^{\frac{1}{2}(\frac{N_\K}{2}+z)}_-(\xi, \z)}\right].\end{align*}
%\end{align*}
 Using the reflection principle of gamma function $\Gm(-z)\Gm(1+z)= \frac{-\pi}{\sin \pi z},$ we get   
\begin{align*}(c^2+P_1\widehat{)^z_+}(\xi, \z)&= \frac{C_{n, d, \K} 2^{z+1}c^{\frac{N_\K}{2}+z}\Gm(1+z)}{\sqrt{|D|}\pi}\left[- \sin\Big(\frac{N_\K}{2}+z\Big)\frac{K_{\frac{N_\K}{2}+z}(c{Q_1}^\frac{1}{2}_+(\xi, \z))}{{Q_1}^{\frac{1}{2}(\frac{N_\K}{2}+z)}_+(\xi, \z)}\right. \\&\qquad \qquad \qquad \qquad \qquad \qquad \qquad \qquad \quad + \left. \frac{\pi}{2} \frac{J_{\frac{N_\K}{2}+z}(c{Q_1}^\frac{1}{2}_-(\xi, \z))}{{Q_1}^{\frac{1}{2}(\frac{N_\K}{2}+z)}_-(\xi, \z)}\right].\end{align*}
Here $\sqrt{|D|}= 1$ and 
$${Q_1}_-^\la(\xi, \z)=\begin{cases}0,~~ &\mbox{if}~~ Q_1(\xi, \z)\geq 0,\\
|Q_1(\xi, \z)|^\la,   &\mbox{if}~~ Q_1(\xi, \z)< 0.\end{cases}$$ 
Now  taking $c= 1$, we get that 
%$$\frac{(1+P_1\widehat{)^z_+}(\xi, \z)}{\Gm(1+z)}=  C_{n, d, \K} 2^{z}\frac{J_{\frac{N_\K}{2}+z}({(|\xi|^2+|\z|^2)}^\frac{1}{2})}{{(|\xi|^2+|\z|^2)}^{\frac{1}{2}(\frac{N_\K}{2}+z)}},$$
% Hence,
\begin{equation}\label{fdgf4}\frac{(1-|x|^2-|y|^2\widehat{)^z_+}(\xi, \z)}{\Gm(1+z)}= \frac{C_{n, d, \K} 2^{z}}{\pi}\frac{J_{\frac{N_\K}{2}+z}({(|\xi|^2+|\z|^2)}^\frac{1}{2})}{{(|\xi|^2+|\z|^2)}^{\frac{1}{2}(\frac{N_\K}{2}+z)}}.\end{equation}
Since    $\frac{J_\la(u)}{u^\la}$ is bounded as $u\to 0,$ for any $\la,$ it is clear that the above function is bounded if and only if $\Re(z)>-\frac{N_\K+1}{2}.$
%-------------------------------------------------------------------------------------------------------------------------------------------------------------

\subsection{Fourier-Dunkl transform of $(1-|x|^2+|y|^2)_+^z$}

Let $\ap_j=-1,~~\bt=1$ for all $j=1, \ldots, n$, $P_2=-|x|^2+|y|^2$ and $Q_2(\xi, \z)= -|\xi|^2+|\z|^2.$ Then $\sqrt{D}= \sqrt{|D|}e^{\frac{1}{2} ni\pi}$ and the formula \eqref{fdgf2} gives
\begin{equation}\label{fd13}
(c^2+P_2+i0\widehat{)^z}(\xi, \z)= \frac{C_{n, d, \K} 2^{z+1}c^{\frac{N_\K}{2}+z}}{\sqrt{|D|}e^{\frac{1}{2} ni\pi}\Gm(-z)}\frac{K_{\frac{N_\K}{2}+z}(c(Q_2(\xi, \z)-i0)^\frac{1}{2})}{(Q_2(\xi, \z)-i0)^{\frac{1}{2}\big(\frac{N_\K}{2}+z\big)}}.
\end{equation}
Similarly,   \eqref{fdgf3} gives
\begin{equation}\label{fd14}
(c^2+P_2-i0\widehat{)^z}(\xi, \z)=  \frac{C_{n, d, \K} 2^{z+1}c^{\frac{N_\K}{2}+z}}{\sqrt{|D|}e^{-\frac{1}{2} ni\pi}\Gm(-z)}\frac{K_{\frac{N_\K}{2}+z}(c(Q_2(\xi, \z)+i0)^\frac{1}{2})}{(Q_2(\xi, \z)+i0)^{\frac{1}{2}\big(\frac{N_\K}{2}+z\big)}}.
\end{equation}
%where $\sqrt{|\D|}e^{-\frac{1}{2}N_\K i\pi}=e^{-\frac{1}{2}N_\K i\pi}.$ 
The above equations \eqref{fd13} and \eqref{fd14} can be expressed in terms of ${Q_2}_+(\xi, \z)$ and ${Q_2}_-(\xi, \z)$ by using the identities \eqref{gi}, \eqref{kpe}, and \eqref{kpe2}, i.e., 

\begin{multline}\label{ffdt3}
(c^2+P_2+i0\widehat{)^z}(\xi, \z)=\frac{C_{n, d, \K} 2^{z+1}c^{\frac{N_\K}{2}+z}}{\sqrt{|D|}e^{\frac{1}{2} ni\pi}\Gm(-z)}  \left[\frac{K_{\frac{N_\K}{2}+z}(c{Q_2}^\frac{1}{2}_+(\xi, \z))}{{Q_2}^{\frac{1}{2}(\frac{N_\K}{2}+z)}_+(\xi, \z)}\right.\\ + \left.\frac{\pi}{2\sin (\frac{N_\K}{2}+z)\pi}\left\{e^{i(\frac{N_\K}{2}+z)\pi}\frac{J_{-(\frac{N_\K}{2}+z)}(c{Q_2}^\frac{1}{2}_-(\xi, \z))}{{Q_2}^{\frac{1}{2}(\frac{N_\K}{2}+z)}_-(\xi, \z)}-\frac{J_{(\frac{N_\K}{2}+z)}(c{Q_2}^\frac{1}{2}_-(\xi, \z))}{{Q_2}^{\frac{1}{2}(\frac{N_\K}{2}+z)}_-(\xi, \z)}\right\}\right],
\end{multline}
and
\begin{multline}\label{ffdt4}
(c^2+P_2-i0\widehat{)^z}(\xi, \z)=\frac{C_{n, d, \K} 2^{z+1}c^{\frac{N_\K}{2}+z}}{\sqrt{|D|}e^{-\frac{1}{2} ni\pi}\Gm(-z)}   \times\left[\frac{K_{\frac{N_\K}{2}+z}(c{Q_2}^\frac{1}{2}_+(\xi, \z))}{{Q_2}^{\frac{1}{2}(\frac{N_\K}{2}+z)}_+(\xi, \z)}\right.\\+\left.\frac{\pi}{2\sin (\frac{N_\K}{2}+z)\pi}\left\{e^{-i(\frac{N_\K}{2}+z)\pi}\frac{J_{-(\frac{N_\K}{2}+z)}(c{Q_2}^\frac{1}{2}_-(\xi, \z))}{{Q_2}^{\frac{1}{2}(\frac{N_\K}{2}+z)}_-(\xi, \z)}-\frac{J_{(\frac{N_\K}{2}+z)}(c{Q_2}^\frac{1}{2}_-(\xi, \z))}{{Q_2}^{\frac{1}{2}(\frac{N_\K}{2}+z)}_-(\xi, \z)}\right\}\right].
\end{multline}
Inserting equations \eqref{ffdt3} and \eqref{ffdt4} in the equation analogous to \eqref{lcfdt}, we obtain
\begin{multline*}
\frac{(c^2+P_2\widehat{)^z_+}(\xi, \z)}{\Gm(1+z)}= \frac{C_{n, d, \K} 2^{z+1}c^{\frac{N_\K}{2}+z}}{\pi}\left[- \sin \Big(z+\frac{n}{2}\Big)\pi \frac{K_{\frac{N_\K}{2}+z}(c{Q_2}^\frac{1}{2}_+(\xi, \z))}{{Q_2}^{\frac{1}{2}(\frac{N_\K}{2}+z)}_+(\xi, \z)}\right. \\+\left. \frac{\pi}{2\sin (\frac{N_\K}{2}+z)\pi} \left(\sin \Big(\frac{d+2\gk}{2}\Big)\pi \frac{J_{-(\frac{N_\K}{2}+z)}(c{Q_2}^\frac{1}{2}_-(\xi, \z))}{{Q_2}^{\frac{1}{2}(\frac{N_\K}{2}+z)}_-(\xi, \z)}+ \sin \Big(z+\frac{n}{2}\Big)\pi \frac{J_{(\frac{N_\K}{2}+z)}(c{Q_2}^\frac{1}{2}_-(\xi, \z))}{{Q_2}^{\frac{1}{2}(\frac{N_\K}{2}+z)}_-(\xi, \z)}\right)\right].
\end{multline*}
The above computation for $c=1,$ can be written as  
\begin{multline}\label{FDTGFH1}
\frac{(1-|x|^2+|y|^2\widehat{)^z_+}(\xi, \z)}{\Gm(1+z)}= \frac{C_{n, d, \K} 2^{z+1}}{\pi}\left[- \sin \Big(z+\frac{n}{2}\Big)\pi \frac{K_{\frac{N_\K}{2}+z}(-|\xi|^2+|\z|^2)^\frac{1}{2}_+)}{{(-|\xi|^2+|\z|^2)}^{\frac{1}{2}(\frac{N_\K}{2}+z)}_+}+ \frac{\pi}{2\sin (\frac{N_\K}{2}+z)\pi}\right.\\ \times\left. \left(\sin \Big(\frac{d+2\gk}{2}\Big)\pi \frac{J_{-(\frac{N_\K}{2}+z)}\Big((-|\xi|^2+|\z|^2)^\frac{1}{2}_-\Big)}{{(-|\xi|^2+|\z|^2)}^{\frac{1}{2}(\frac{N_\K}{2}+z)}_-}+ \sin \Big(z+\frac{n}{2}\Big)\pi \frac{J_{(\frac{N_\K}{2}+z)}(-|\xi|^2+|\z|^2)^\frac{1}{2}_-)}{{(-|\xi|^2+|\z|^2)}^{\frac{1}{2}(\frac{N_\K}{2}+z)}_-}\right)\right].
\end{multline}

%An immediate application of the above calculation is the following one.
%\begin{lem}\label{fdp}
%	If $P(x, y)=|x|^2-|y|^2$ for $(x, y)\in\R^n\times\R^d$ then
%	\begin{align}\label{fdpz}\nonumber
%		\widehat{P^z_+}(\xi, \z)&=c_\K^{-1}\pi^{-\frac{n}{2}}2^{2z+\frac{d}{2}+\gk}\Gamma\Big(z+\frac{N_k}{2}\Big)\Gamma(1+z)\\ &\qquad \times \left(e^{-\frac{1}{2}i\pi(d+2\gk)}e^{-i\pi z}(Q-i0)^{-z-\frac{N_k}{2}}-e^{\frac{1}{2}i\pi(d+2\gk)}e^{i\pi z}(Q+i0)^{-z-\frac{N_k}{2}}\right).
%	\end{align}
%\end{lem}

\section{Proof of main results}\label{pmr}
This section is devoted to provide    proofs of our main results.
\begin{proof}[\bf Proof of Theorem \ref{FDRTP} and \ref{REOp}]
Let $G_z$ denotes the generalized function corresponding (see Section \ref{gtrt}) to the paraboloid-surface $S_1$, defined by
\be \label{eq:gfforparaboloid} G_z(x, y)=(x_n-|x'|^2+|y|^2)_+^z, \quad (x, y)\in \R^n \times \R^d,\ee where $x'=(x_1, \ldots, x_{n-1})$. In view of \eqref{DFgp1}, $\widehat{G_z}(\xi, \z)$ is bounded in $(\xi, \z)$ if and only if $\Re(z)= - \frac{(d+n+2\gk+1)}{2}$. Thus the proof of the theorems will follow from Theorem \ref{C}.
\end{proof}

%\begin{theorem}[Restriction estimates for orthonormal functions-cone case]\label{reoc}
%Let $S$ be the surface given by $S=\{(x, y)\in\R^n\times\R^d: |x|^2-|y|^2=0\},~~~ r=\frac{N_k}{N_k-2},~~ \bt=\frac{N_k}{N_k-1}$ and $N_k\geq 2.$ Then for any (possible infinite) orthonormal system $f_j$ in $L^2(S, h^2_\K(y)d\sigma (x, y))$ and for any $(n_j)\subset \C$ we have
%\begin{equation}\label{rec}
%\left\|\sum_jn_j|\ecs f_j|^2\right\|_{L^r_\K(\R^n\times\R^d)}\leq C \left(\sum_j|n_j|^\bt\right)^\frac{1}{\bt}
%\end{equation}
%\end{theorem}
\begin{proof}[{\bf Proof of Theorem \ref{FDRTS2} and \ref{REOs}}]
For the sphere $S_2$, we consider the corresponding generalized function defined  by 
$$G_z(x, y)=\frac{1}{\Gm(1+z)}(1-|x|^2-|y|^2)^z_+.$$ In view of  \eqref{fdgf4},  $\widehat{G_z}(\xi, \z)$ is bounded in $(\xi, \z)$ if and only if $\Re(z)>-\frac{N_\K+1}{2}.$ Thus the proof of the theorems will follow from Theorem \ref{C}.
\end{proof}

\begin{proof}[{\bf Proof of Theorem \ref{FDRTh} and \ref{REOh}}]
For the two sheet hyperboloid $S_3$, we let 
$$G_z(x, y)=\frac{h(z)}{\Gm(1+z)}(1-|x|^2+|y|^2)^z_+,$$where $h(z)$ will be specified later.

%$$\mbox{let}\quad G_z(x, y)=h(z)\Gm(z+1)^{-1}(1+P_2)^z_+,$$
%where $h(z)$ will be specified later. Then the Fourier-Dunkl transform of $G_z(x, y)$ is 
%$$\widehat G_z(\xi, \z)= h(z)\Gm(z+1)^{-1}(1+P_2\widehat{)^z_+}(\xi, \z).$$

Case I.
Let $n=1=d,$ and $\gk=0$ so that $N_\K= 2.$ Choose $h(z)= (z+1)^{-1} \sin \pi(z+1).$  Since,  for $z=-1$ the term $  \frac{K_{\frac{N_\K}{2}+z}({Q_2}^\frac{1}{2}_+)}{{Q_2}^{\frac{1}{2}(\frac{N_\K}{2}+z)}_+}$ has a pole, in view of \eqref{FDTGFH1},  $\widehat G_z$ is bounded only for $-1>\Re (z)\geq -\frac{3}{2}.$ 

Case II.
Let $n\geq 1,~~d\geq 1$ and $N_\K>2,$ we choose  the function $h$ as
$$h(z)=\begin{cases}
	 \Big(z+\frac{N_\K}{2}\Big)(z+1)^{-1} \sin \Big(z+\frac{N_\K}{2}\Big)\pi,& \text{if $N_\K$ is an even integer},\\ \Big(z+\frac{N_\K}{2}\Big) \sin \Big(z+\frac{N_\K}{2}\Big)\pi,& \text{Otherwise}. 
\end{cases}$$
 Note that $h(-1)\neq 0$. Thus  $G_{-1}$ is a non-zero multiple of $d\mu_r.$ We prove that $\widehat{G_z}$ is bounded if $-\frac{N_\K+1}{2}\leq\Re (z)\leq -\frac{N_\K}{2}.$ The multiplication with  $h(z)$  cancels the poles of $\displaystyle \frac{1}{\sin \big(z+\frac{N_\K}{2}\big)\pi}$ in this region.   

Let $\la= -z-\frac{N_\K}{2}$. Then in this case, the region of $z$ will become 
\begin{equation}\label{rexp}
0\leq\Re (\la)\leq\frac{1}{2}.
\end{equation}
On this region of $\la,$ it is enough to prove the boundedness of $u^\la J_{\la}(u),~~u^\la J_{-\la}(u),$ and $\la u^\la K_{\la}(u).$ We note that $K_\la=K_{-\la}$.  Also, $u^\la J_{\la}(u),$ and $u^\la J_{-\la}(u)$ are bounded on the region given by \eqref{rexp} because of the well known following Bessel functions estimates
$$J_\la(u)= \begin{cases} O(u^\la);~~~&u\to 0\\O(u^{-1/2});~~~&u\to\infty\end{cases}$$
for the term $\la u^\la K_\la(u)$ we use the integral formula $\displaystyle K_\la(u)= \int_0^\infty \cosh \la t\; e^{-u \cosh t} dt$ to obtain exponential decay as $u\to\infty$ \cite[p. 259]{JV},  and the power series expansion \cite[p. 270]{JV},  
$$\la u^\la K_\la(u)= \frac{1}{2}\sum_{k=0}^\infty \frac{(-1)^k}{k!}[\la \Gm(-k-\la)(u/2)^{2\la}+ \la \Gm(\la-k)] (u/2)^{2k},$$
to obtain the boundedness for small values of $u.$ Notice that the factor $\la$ cancels the poles of the $\Gm$-function at $\la=0.$
In view of the above analysis and Theorem \ref{C}, proof of the theorems will follow.

\end{proof}

%------------------------------------------------------------------------------------------------------------------------------------------------------------

\section{Strichartz inequalities}\label{sec5}
This section is devoted to study Strichartz inequality for a system of orthonormal functions associated  with  Dunkl Laplacian propagator $e^{it\Dk}$  and Klein-Gordon propagator $e^{it\sqrt{1-\Dk}}$.  We begin with       Dunkl Laplacian  case.
\subsection{Dunkl Laplacian case} Strichartz inequalities are important applications of the restriction estimates of quadratic surfaces, which are useful tools to study non-linear Schr\"odinger, wave, and Klein-Gordon equations.

Let $S$ be the paraboloid given by $$S=\{(\om, \z)\in\R\times\R^d: \om= -|\z|^2\}$$ with the measure $d\mu(\om, \z)=d\z$.  Then   for all $f\in L^1(S, d\mu)$ and for all $(t, y)\in\R\times\R^d,$ we have  
\begin{align*}
\ecs f(t, y)&=\frac{1}{c_\K\sqrt{2\pi}}\int_Sf(\om, \z)e^{it\cdot\om}E_\K(i\z, y)h^2_\K(\z)d\mu(\om, \z)\\
&=\frac{1}{c_\K\sqrt{2\pi}}\int_{\R^d}f(-|\z|^2, \z)e^{-it\cdot|\z|^2}E_\K(i\z, y)h^2_\K(\z)d\z.
\end{align*}
In particular, by choosing $f(\om, \z)= \mathcal{F}_\K\phi(\z),$ for some $\phi: \R^d\to \C,$ we get 
\begin{align}
\ecs f(t, y)&=\frac{1}{c_\K\sqrt{2\pi}}\int_{\R^d}\mathcal{F}_\K\phi(\z)e^{-it\cdot|\z|^2}E_\K(i\z, y)h^2_\K(\z)d\z\nonumber \\
&=\frac{1}{\sqrt{2\pi}}(e^{it\Dk}\phi)(y), \label{restrel1}
\end{align}
for all $(t, y)\in\R\times\R^d.$

Now  consider the Cauchy problem for the free Schr\"odinger equation associated with the Dunkl Laplacian operator $\Delta_{k}$, namely
\begin{align}\label{100}
	\left\{  \begin{array}{ll}i \partial_{t} u(t, x)+ \Delta_{\K} u(t, x)=0, & (t, x) \in \mathbb{R} \times  \mathbb{R}^d, \\ u(0,x)=\phi(x),& x \in   \mathbb{R}^d.  \end{array}\right.
\end{align}
For $\phi\in L^2(\R^d, h^2_\K(y)dy)$,  $u(t, x)=e^{i t \Delta_\K} \phi(x)$ is the solution of the above system.   Using Theorem \ref{FDRTP} and duality argument we get the following Strichartz estimates
$$\|e^{it\Dk}\phi\|_{L^{p'}_\K(\R\times\R^d)}\leq C \|\phi\|_{\K, 2},$$
for all $\phi\in L^2(\R^d, h^2_\K(y)dy),$ where $p'= 2+\frac{4}{d+2\gk}.$

\begin{theorem}[Strichartz estimates for orthonormal functions-diagonal case]\label{SEOpd}
Let $d\geq 1,$ then for any (possibly infinite) orthonormal system $(\phi_j)$ in $L^2(\R^d, h^2_\K(y) dy)$ and for any $(n_j)\subset \C$, we have 
$$\left\|\sum_j n_j|e^{it\Dk}\phi_j|^2\right\|_{L^\frac{p'}{2}_\K(\R\times \R^d)}\leq C \Big(\sum_j|n_j|^\bt\Big)^\frac{1}{\bt},$$
where $\bt= \frac{d+2\gk+2}{d+2\gk+1}$ and $p'= 2+\frac{4}{d+2\gk}$ with $C>0$ independent of $(n_j)$ and $(\phi_j).$
\end{theorem}
\begin{proof}
Let $f_j(\om, \z)=\mathcal{F}_\K(\phi_j)(\z)$ for all $\z\in \R^d$ and $\om=-|\z|^2$. Since $(\phi_j)$ is an orthonormal system in $L^2(\R^d, h^2_\K(y) dy)$, $(f_j)$ is an orthonormal system in $L^2(S, h_\K^2(\z)d\sigma(\om, \z))$. In view of \eqref{restrel1} and Theorem \ref{REOp} with $n=1$, we get the required inequality.
\end{proof}

\begin{theorem}[Strichartz estimates for orthonormal functions-general case]\label{SEOpg}
Let $d\geq 1$ and $p, q\geq 1$ such that $\frac{2}{p}+\frac{d+2\gk}{q}=d+2\gk,~~ 1\leq q<1+\frac{2}{d+2\gk-1}.$ Then for any (possible infinite) orthonormal system $\phi_j$ in $L^2(\R^d, h^2_\K(y) dy)$ and for any $(n_j)\subset \C$, we have 
$$\Big\|\sum_j n_j|e^{it\Dk}\phi_j|^2\Big\|_{L^p_t(\R, L^q_y(\R^d, h^2_\K(y) dy))}\leq C\Big(\sum_j |n_j|^\frac{2q}{q+1}\Big)^\frac{q+1}{2q},$$
with $C>0$ independent of $(n_j)$ and $(\phi_j).$ 
\end{theorem}
\begin{proof}
Proof will follow from the Schatten bound with mixed norm (Theorem \ref{sbmt}) along with duality principle in mixed norm space (Lemma \ref{dtm}). 
\end{proof}

\begin{lem}[Duality principle in mixed norm space]\label{dtm}
Let $p\geq 1,~q\geq 1$ and $\ap\geq 1.$ Let $A$ be a bounded linear operator from a Hilbert space $\Hc$ to $L_t^{2p}(\R, L_y^{2q}(\R^d, h^2_\K(y) dy)).$ Then the following statements are equivalent.

$(1)$ There is a constant $C > 0$ such that
	\begin{equation}\label{DI}
		\|WAA^*\overline{W}\|_{\s^\ap\left(L^2_{\K}(\R\times \R^d)\right)} \leq C \|W\|^2_{L^{2p'}_t(\R, L^{2q'}_y(\R^d, h^2_\K(y) dy))},
	\end{equation}
	for all $W\in {L^{2p'}_t(\R, L^{2q'}_y(\R^d, h^2_\K(y) dy))},$ where the function $W$ is interpreted as an operator which acts by multiplication.\\
	
	$(2)$ For any orthonormal system $(f_j)_{j\in J}$ in $\Hc$ and any sequence $(n_j)_{j\in J}\subset \C,$ we have
	
	\begin{equation}\label{DI2}
		\left\|\sum_{j\in J}n_j\left|Af_j\right|^2\right\|_{L^{p}_t(\R, L^{q}_y(\R^d, h^2_\K(y) dy))}\leq C' \left(\sum_{j\in J}|n_j|^{\ap'}\right)^{1/\ap'},
	\end{equation}
where $C'$ is a constant. Moreover, the values of the optimal constants $C$ and $C'$ coincide.
\end{lem}
\begin{proof}
Proof of this lemma is similar to Lemma \ref{B}.
\end{proof}

\begin{theorem}[Schatten bound with mixed norms]\label{sbmt}
Let $d\geq 1$ and $S$ be the paraboloid given by $$S:= \{(\om, \z)\in\R\times\R^d: \om=-|\z|^2\}.$$ Then, for all exponents $p, q\geq 1$ satisfying the relations $$\frac{2}{p}+\frac{d+2\gk}{q}=d+2\gk$$ and $$1+\frac{2}{d+2\gk}<q<\frac{d+2\gk+1}{d+2\gk-1},$$ we have 
$$\|W_1T_SW_2\|_{\s^{2q'}(L^2_\K(\R\times\R^d))}\leq C\|W_1\|_{L^{2p'}_t(\R, L^{2q'}_y(\R^d, h^2_\K(y) dy))}\|W_2\|_{L^{2p'}_t(\R, L^{2q'}_y(\R^d, h^2_\K(y) dy))},$$
with $C>0$ independent of $W_1, W_2.$
\end{theorem}
\begin{proof}
Consider the generalized function $G_z$ for the given surface (see Section \ref{gtrt}); $$G_z(\om, \z)= \frac{1}{\Gm(z+1)}(\om+|\z|^2)^z_+,\;\forall (\om, \z)\in\R\times\R^d,$$
which ensures that the Fourier-Dunkl multiplication operator $(T_z)$ with $G_{-1}$ coincides with the operator $T_S.$ In \cite[Section 3]{JPSH}, it is observed that  
$$\|T_{is}\|_{L^2_\K(\R^{d+1})\rightarrow L^2_\K(\R^{d+1})}=\|G_{is}\|_{L^\infty(\R^{d+1})}\leq \left|\frac{1}{\Gm(1+is)}\right|\leq C e^{\pi|s|/2}.$$
In view of \eqref{DFgp1}, we have 
\begin{align}
\check{G}(t, y)= & C_\K i e^{iz\frac{\pi}{2}} e^{-\frac{i\pi}{4}(d+2\gk)\sgn t } e^{\frac{i}{4t}|y|^2} |t|^{-\frac{(d+2\gk)}{2}}(t+i0)^{-z-1},\; \forall(t, y)\in\R\times\R^d,
\end{align}
where $C_\K=\frac{1}{c_\K\sqrt{\pi}2^\frac{d+2\gk-n+2}{2}}.$

For diagonal case (Theorem \ref{SEOpd}), we just used the fact that $\check{G}(t, y)$ belongs to $L^\infty_\K(\R\times\R^d),$ but for this theorem (general case) we need better estimates of $\check{G}(t, y).$  To do so, recall that the distribution $(t+i0)^\la$ on $\R$ satisfies the identity  (see \cite[Ch. 1, Sec. 3.6]{GS})
$$(t+i0)^\la=t^\la_+ + e^{i\pi\la}t^\la_-,$$
for $\Re (\la)>-1$, where $t^\la_{\pm}$ are the distributions given by the $L^1_{loc}$-functions 
$$t^\la_+=\begin{cases} t^\la, &\mbox{for}~ t>0,\\ 0, &\mbox{for}~ t\leq 0,\end{cases}\qquad\text{and }\qquad t^\la_-=\begin{cases} 0 ,&\mbox{for}~ t\geq 0,\\ (-t)^\la, &\mbox{for}~ t< 0.\end{cases}$$
In particular, the distribution $(t+i0)^\la$ is also given by a $L^1_{loc}$-function, and we deduce the bound 
$$|(t+i0)^\la|\leq \max \Big\{1, e^{-\pi\Im \la}\Big\}|t|^{\Re( \la)}, \; \forall t\in\R,$$ valid for all $\Re( \la)>-1.$ In our context, we have $\la= -z-1$ with $z=-\la_0+is,$ so that $\Re (\la)= \la_0-1>0.$  For all $s\in\R$ and for all $\la_0>1,$ we   deduce that
\begin{equation}\label{ings}
|\check{G}_{-\la+is}(t, y)|\leq C e^{\pi |s|/2}|t|^{\la_0-1-(d+2\gk)/2}, 
\end{equation}
for all $(t, y)\in\R\times\R^d.$ Using an application of Hardy-Littlewood-Sobolev inequality (see \cite{BW}) along with \eqref {ings} yields  
\begin{align*}
&\|W^{\la_0-is}_1T_{-\la_0+is}W^{\la_0-is}_2\|^2_{\s^2(L^2_\K(\R\times\R^d))}\\
&=\int_{\R^{2(d+1)}}|W_1(t, y)|^{2\la_0}|\tau_y \check{G}_{-\la_0+is}(t-t', y')|^2 |W_2(t', y')|^{2\la_0}h^2_\K(y) h^2_\K(y')  dy dy' dt dt'\\
&\leq C\int_{\R^{2(d+1)}}|W_1(t, y)|^{2\la_0}\|\check{G}_{-\la_0+is}(t-t', y')\|^2_{\K, \infty} |W_2(t', y')|^{2\la_0}h^2_\K(y) h^2_\K(y')  dy dy' dt dt'\\
&\leq C e^{\pi|s|}\int_{\R}\int_{\R} \frac{\|W_1(t, \cdot)\|^{2\la_0}_{\K, 2\la_0} \|W_2(t', \cdot)|^{2\la_0}_{\K, 2\la_0}}{|t-t'|^{d+2\gk+2-2\la_0}} dt dt'\\
&\leq C e^{\pi|s|}\|W_1\|^{2\la_0}_{L^{\frac{4\la_0}{2\la_0-d-2\gk}}_t(\R, L^{2\la_0}_y(\R^d, h^2_\K(y)dy))} \|W_2\|^{2\la_0}_{L^{\frac{4\la_0}{2\la_0-d-2\gk}}_t(\R, L^{2\la_0}_y(\R^d, h^2_\K(y)dy))},
\end{align*}
provided we have $0\leq d+2\gm_k+2-2\la_0 < 1,$ i.e.,  $\frac{(d+2\gk+1)}{2} < \la_0 \leq \frac{d+2\gk+2}{2}.$ By Theorem 2.9 of \cite{S}, we have
$$\|W_1T_{-1}W_2\|_{\s^{2\la_0}(L^2_\K(\R\times\R^d))} \leq C \|W_1\|_{L^{\frac{4\la_0}{2\la_0-d-2\gk}}_t(\R, L^{2\la_0}_y(\R^d, h^2_\K(y)dy))} \|W_2\|_{L^{\frac{4\la_0}{2\la_0-d-2\gk}}_t(\R, L^{2\la_0}_y(\R^d, h^2_\K(y)dy))}, $$ for $\frac{d+2\gk+1}{2}<\la_0\leq \frac{d+2\gk+2}{2}.$ Since $ \frac{d+2\gk+2}{d+2\gk}<q<\frac{d+2\gk+1}{d+2\gk-1}$, we have $\frac{d+2\gk+1}{2}<q'<\frac{d+2\gk+2}{2}$. Thus the theorem is proved by choosing $\la_0=q'$.

\end{proof}

%------------------------------------------------------------------------------------------------------------------------------------------------------------

\subsection{Klein-Gordon case}
In this subsection, we prove Strichartz inequalities for the Klein-Gordon propagator $e^{it\sqrt{1-\Dk}}$ associated with Dunkl Laplacian using Theorem \ref{FDRTh} and Theorem \ref{REOh}.  

Let $S=\{(\om, \z)\in\R\times\R^d: \om^2=1+|\z|^2\}$ with the measure $d\mu(\om, \z)=\frac{d\z}{2\sqrt{1+|\z|^2}}.$ Then  for  $f\in L^1(S, h^2_\K(\z)d\mu)$ and for all $(t, y)\in\R\times\R^d,$    we obtain
\begin{align*}
\ecs f(t, y)=&\frac{1}{c_\K\sqrt{2\pi}}\int_Sf(\om, \z)e^{it\cdot\om}E_\K(i\z, y)h^2_\K(\z)d\mu(\om, \z)\\
=&\frac{1}{c_\K\sqrt{2\pi}}\int_{\R^d}f(\sqrt{1+|\z|^2}, \z)e^{it\cdot\sqrt{1+|\z|^2}}E_\K(i\z, y)h^2_\K(\z)\frac{d\z}{2(\sqrt{1+|\z|^2})}\\ &\quad + \frac{1}{c_\K\sqrt{2\pi}}\int_{\R^d}f(-\sqrt{1+|\z|^2}, \z)e^{-it\cdot\sqrt{1+|\z|^2}}E_\K(i\z, y)h^2_\K(\z)\frac{d\z}{2(\sqrt{1+|\z|^2})}.
\end{align*}
In particular, if we choose $f(\om, \z)= 2\hspace{0.2cm}\mathds{1}_{(\om>0)} \sqrt{1+|\z|^2} \mathcal{F}_\K(\phi)(\z),$ then we obtain
$$\ecs f(t, y)=\frac{1}{\sqrt{2\pi}}e^{it\sqrt{1-\Dk}}\phi(y),\hspace{0.5cm}\forall (t, y)\in(\R\times\R^d).$$ 
Now  consider the Klein-Gordon equation
\begin{align}\label{100}
	\left\{  \begin{array}{ll} \partial_{tt} u(t, x)+u(t, x)= \Delta_{k} u(t, x), & (t, x) \in \mathbb{R} \times  \mathbb{R}^d, \\ u(0,x)=u_0(x),& x \in   \mathbb{R}^d,\\
		\partial_{t}  u(0,x)=u_1(x),& x \in   \mathbb{R}^d.  \end{array}\right.
\end{align}
Then the solution of the above system can   be written as  $u=u_{+}+u_{-}$, where $u_{+}$ and $u_{-}$ are given by
$$
u_{\pm}(t, x)=e^{\pm i t \sqrt{1-\Dk}} \phi_{\pm}(x),
$$
and $\phi_{+}$ and $\phi_{-} $ satisfy the following
$$\begin{cases}
	u_0=\phi_{+}+\phi_{-},\\
	u_1= i \sqrt{1-\Dk}\left(\phi_{+}-\phi_{-}\right).
\end{cases} $$  
Consequenly, the Strichartz estimates for the   Klein-Gordon equation are usually given by those for the one-sided propagator $e^{\pm i t \sqrt{1-\Dk}}$.  
Due to Theorem \ref{FDRTh}, $\ecs$ is bounded from $L^2(S, h^2_\K(\z)d\mu)$ to $L^{p'}_\K(\R\times\R^d)$ with $2+\frac{4}{d+2\gk}\leq p'\leq 2+\frac{4}{d+2\gk-1},$ when $d+2\gk>1$ and $6\leq p'\leq \infty$ when $d+2\gk=1.$ Therefore we get that 
\begin{equation}\label{equ:StKG}
\|e^{it\sqrt{1-\Dk}}\phi\|_{L^{p'}_\K(\R\times\R^d)}\leq C\|\phi\|_{H_\K^{1/2}(\R^d)},
\end{equation} 
where $H_\K^{1/2}(\R^d)$ is the Dunkl Sobolev space defined in the preliminary section.  Using Theorem  \ref{REOh} we also get the following corresponding version of the estimate \eqref{equ:StKG} for orthonormal functions. 
\begin{theorem}[Strichartz estimates for orthonormal functions- Klein-Gordon case]\label{SEOKG}   
Assume $d\geq 1.$ Let $1+\frac{2}{d+2\gk}\leq r\leq 1+\frac{2}{d+2\gk-1}$ if $d+2\gk>1$ and $3\leq r<\infty$ if $d+2\gk=1.$ For any (possibly infinite) orthonormal system $(\phi_j)$ in $H_\K^{1/2}(\R^d)$ and for any $(n_j)\subset \C,$ we have 
$$\Big\|\sum_j n_j|e^{it\sqrt{1-\Dk}}\phi_j|^2\Big\|_{L^r(\R\times\R^d)}\leq C \Big(\sum_j|n_j|^\frac{2r}{r+1}\Big)^\frac{r+1}{2r},$$
with $C>0$ independent of $(n_j)$ and $(\phi_j).$
\end{theorem}
\begin{proof}
For the given orthonormal system $(\phi_j)$ in $H_\K^{1/2}(\R^d)$, we define $$f_j(\om, \z)=2~\mathds{1}_{(\om>0)} \sqrt{1+|\z|^2} \mathcal{F}_\K(\phi_j)(\z)$$ for $(\om, \z)\in S$. Then $(f_j)$ is an orthonormal system in $L^2(S, h_\K^2(y) d\sigma(\om, \z))$, as we have seen in the beginning of this subsection. We then apply Theorem  \ref{REOh} to this system for $n=1$.
\end{proof}

\begin{cor}
Let $d\geq 1,$ suppose $1+\frac{2}{d+2\gk}\leq r_0\leq 1+\frac{2}{d+2\gk-1},$ for $1\leq r\leq r_0,~~q\geq 1$ and $s\geq 0$ such that $\frac{1}{q}=\frac{1}{1-r_0}\big(1-\frac{1}{r}\big)$ and $s=\frac{r_0}{r_0-1}\big(\frac{1}{2}-\frac{1}{2r}\big).$ Then for all families of orthonormal functions $(\phi_j)$ in $H^s_\K(\R^d),$ we have that 
\begin{equation}\label{kgec}
\Big\|\sum_j n_j|e^{it\sqrt{1-\Dk}}\phi_j|^2\Big\|_{L^q_t(\R, L^r_y(\R^d, h^2_\K(y) dy)}\leq C \Big(\sum_j|n_j|^\bt\Big)^\frac{1}{\bt},
\end{equation}
where $\bt= \frac{2r}{r+1}.$
\end{cor}
\begin{proof}
In view of Theorem \ref{SEOKG}, we have the estimate \eqref{kgec} for the case $(q, r, s, \bt)= \big(r_0, r_0, \frac{1}{2}, \frac{2r_0}{r_0+1}\big),$ Since $e^{it\sqrt{1-\Dk}}$ is an unitary operator on $L^2(\R^d, h^2_\K(y) dy)$ and we have the estimate \eqref{kgec} for the case $(q, r, s, \bt)= (\infty, 1, 0, 1).$ By using the interpolation between these two points we get the required estimate.    
\end{proof}
%------------------------------------------------------------------------------------------------------------------------------------------------------------
\section{Perspective}
In this paper, we considered the Fourier-Dunkl transform and proved Strichartz's restriction theorem for this transform for the paraboloid, sphere, and hyperboloid surfaces. Further,  as an application of the restriction theorem, we proved Strichartz inequality for orthonormal families of initial data. Moreover,  we derive the Strichartz inequality for Schr\"odinger equation associated with Dunkl Laplacian  and Klein-Gordon operator for the family of orthonormal functions.  In the forthcoming paper, we will use the results proved in this preprint for investigate the non-linear Klein-Gordon equation in the Dunkl setting.
%------------------------------------------------------------------------------------------------------------------------------------------------------------

\section*{Acknowledgments}
The authors are deeply indebted to the  professors S. Thangavelu and M. W. Wong
  for their help.

\end{document}